\documentclass[11pt, a4paper]{article}

\usepackage[utf8]{inputenc}
\usepackage[english]{babel}
\usepackage[T1]{fontenc}

\usepackage{amsmath}
\usepackage{amssymb}
\usepackage{amsthm}
\usepackage{mathabx}
\usepackage{tikz}
\usepackage{tikz-cd}
\usepackage{todonotes}
\usepackage{hyperref}
\usepackage{cleveref}
\usepackage{footmisc}
\usepackage{bbm}
\usepackage{afterpage}
\usepackage{mathrsfs}
\usepackage{sseq}
\usepackage{spectralsequences}
\usepackage{extarrows}
\usepackage{enumitem}
\usepackage{stmaryrd}
\usepackage{pdflscape}

\newcommand{\old}{\mathrm{old}}

\newcommand{\ueq}{\overset{\circ}{=}}
\newcommand{\uequiv}{\overset{\circ}{\equiv}}

\newcommand{\1}{\mathbf{1}}
\newcommand{\BP}{\mathrm{BP}}

\DeclareMathOperator{\K}{K}

\DeclareMathOperator{\ko}{ko}

\newcommand{\Sph}{\mathbf{S}}

\DeclareMathOperator{\TMF}{TMF}
\DeclareMathOperator{\tmf}{tmf}

\DeclareMathOperator{\CAlg}{CAlg}

\DeclareMathOperator{\Sp}{Sp}

\DeclareMathOperator{\Syn}{Syn}

\newcommand{\id}{\mathrm{id}}

\newcommand{\E}{\mathbf{E}}

\newcommand{\F}{\mathbf{F}}

\newcommand{\h}{\mathrm{h}}

\newcommand{\J}{\mathrm{J}}
\renewcommand{\j}{\mathrm{j}}

\newcommand{\Z}{\mathbf{Z}}

\newcommand{\al}{\alpha}
\newcommand{\be}{\beta}
\newcommand{\ga}{\gamma}

\usepackage{mathtools}
\DeclarePairedDelimiter{\ceil}{\lceil}{\rceil}

\theoremstyle{theorem}\numberwithin{equation}{section}
\newtheorem{theorem}[equation]{Theorem}
\crefname{theorem}{{th}.\!\!}{{ths}.\!\!}
\Crefname{theorem}{{Th}.\!\!}{{Ths}.\!\!}
\newtheorem{theoremalph}{Theorem}

\crefname{theoremalph}{{th}.\!\!}{{ths}.\!\!}
\Crefname{theoremalph}{{Th}.\!\!}{{Ths}.\!\!}

\Crefname{problem}{{Prb}.\!\!}{{Prbs}.\!\!}
\newtheorem{prop}[equation]{Proposition}
\Crefname{prop}{{Pr}.\!\!}{{Prs}.\!\!}
\newtheorem{lemma}[equation]{Lemma}
\Crefname{lemma}{{Lm}.\!\!}{{Lms}.\!\!}
\newtheorem{cor}[equation]{Corollary}
\Crefname{cor}{{Cor}.\!\!}{{Cors}.\!\!}

\Crefname{conjecture}{{Conj}.\!\!}{{Conjs}.\!\!}

\theoremstyle{definition}\numberwithin{equation}{section}
\newtheorem{mydef}[equation]{Definition}
\Crefname{mydef}{{Df}.\!\!}{{Dfs}.\!\!}

\Crefname{variant}{{Var}.\!\!}{{Vars}.\!\!}

\Crefname{recall}{{Rcl}.\!\!}{{Rcls}.\!\!}

\Crefname{construction}{{Con}.\!\!}{{Cons}.\!\!}

\Crefname{ass}{{As}.\!\!}{{As}.\!\!}
\newtheorem{question}[equation]{Question}
\Crefname{question}{{Q}.\!\!}{{Qs}.\!\!}

\Crefname{notation}{{Nt}.\!\!}{{Nts}.\!\!}

\Crefname{situation}{{St}.\!\!}{{Sts}.\!\!}

\theoremstyle{remark}\numberwithin{equation}{section}

\Crefname{example}{{Ex}.\!\!}{{Exs}.\!\!}

\Crefname{nonexample}{{NonEx}.\!\!}{{NonEx}.\!\!}

\Crefname{claim}{{Clm}.\!\!}{{Clms}.\!\!}
\newtheorem{remark}[equation]{Remark}
\Crefname{remark}{{Rmk}.\!\!}{{Rmks}.\!\!}

\Crefname{idea}{{Id}.\!\!}{{Ids}.\!\!}

\Crefname{warn}{{Warn}.\!\!}{{Warns}.\!\!}

\Crefname{figure}{{Fig.}\!\!}{{Figs.}\!\!}
\Crefname{footnote}{{Fn.}\!\!}{{Fn.}\!\!}

\Crefname{part}{{\textsection}\!\!}{{\textsection}\!\!}
\Crefname{chapter}{{\textsection}\!\!}{{\textsection}\!\!}
\Crefname{section}{{\textsection}\!\!}{{\textsection}\!\!}
\Crefname{subsection}{{\textsection}\!\!}{{\textsection}\!\!}
\Crefname{appendix}{{\textsection}\!\!}{{\textsection}\!\!}

\begin{sseqdata}[name = ANSSm3, Adams grading, classes = fill, tick step = 2,
major tick step = 4,
minor tick step = 1]    

\class[red](72,0)\class[red](76,0)\class[red](80,0)\class[red](84,0)

\class(71,1)\class[blue](75,1)\class[red](75,1)\class[blue](79,1)\class[blue](83,1)\class[red](83,1)

\class[blue](74,2)\class[red](82,2)\class[red](82,2)\class[blue](82,2)\class(86,2)

\class[blue](81,3)\class[blue](81,3)\class(85,3)\class(85,3)\class(85,3)

\class[red](76,4)\class(84,4)\class(84,4)\class(84,4)

\class[blue](75,5)\class[red](79,5)\class[red](79,5)\class[blue](83,5)

\class[red](78,6)\class[blue](78,6)\class[blue](78,6)\class[red](82,6)

\class[red](73,7)\class[blue](77,7)\class[red](81,7)\class[blue](81,7)

\class[blue](72,8)\class[red](76,8)\class[blue](80,8)

\class[blue](75,9)\class[red](75,9)

\class[blue](74,10)\class[red](78,10)\class[blue](82,10)\class[blue](86,10)

\class[red](77,11)\class[blue](77,11)\class(85,11)

\class[blue](76,12)\class[red](80,12)

\class[red](71,13)\class[blue](79,13)

\class[blue](70,14)\class[red](74,14)

\class[blue](73,15)\class[red](81,15)

\class[blue](80,16)\class(84,16)

\class(83,17)

\structline[orange](72,0)(75,1)\structline[orange](76,0)(79,1)\structline[orange](80,0)(83,1)
\structline[orange](82,2,1)(85,3,1)\structline[orange](82,2,2)(85,3,2)\structline[orange](82,2,3)(85,3,3)
\structline[orange](81,3,1)(84,4,1)\structline[orange](81,3,2)(84,4,2)
\structline[orange](76,4)(79,5,1)
\structline[orange](75,5)(78,6,2)\structline[orange](75,5)(78,6,3)\structline[orange](79,5,2)(82,6)
\structline[orange](78,6,1)(81,7,1)\structline[orange](78,6,2)(81,7,2)\structline[orange](78,6,3)(81,7,2)
\structline[orange](73,7)(76,8)\structline[orange](77,7)(80,8)
\structline[orange](72,8)(75,9,1)
\structline[orange](75,9,2)(78,10)
\structline[orange](74,10)(77,11,2)\structline[orange](82,10)(85,11)
\structline[orange](77,11,1)(80,12)
\structline[orange](76,12)(79,13)
\structline[orange](71,13)(74,14)
\structline[orange](70,14)(73,15)
\structline[orange](81,15)(84,16)
\structline[orange](80,16)(83,17)

\end{sseqdata}

\begin{sseqdata}[name = homotopyoftmfpsiat3, Adams grading, tick step = 4,
major tick step = 4,
minor tick step = 1]    

\foreach \i in {1,2,3,4} {
	\class[fill=blue](10*\i, 2*\i)
	}

\class[fill=blue, rectangle](0,0)

\foreach \i in {0,1,2,3} {
	\structline[fill=black](10*\i, 2*\i)(10*(\i+1),2*(\i+1))
	}

\foreach \i in {0,1,2,3} {
\class[fill=blue](3+72*\i,1)
\class[fill=blue](13+72*\i,3)
\class[fill=blue](27+72*\i,1)
\class[fill=blue](37+72*\i,3)
	}

\structline[black](0,0)(3,1)
\structline[black](10,2)(13,3)
\structline[black](3,1)(13,3)
\structline[black](27,1)(30,6)
\structline[black](27,1)(37,3)
\structline[black](37,3)(40,8)

\class[fill=red](7,1)
\class[fill=red, circlen = 2](11,1)
\class[fill=red](15,1)
\class[fill=red](19,1)
\class[fill=red, circlen = 2](23,1)
\class[fill=red](23,1)
\class[fill=red](27,1)
\class[fill=red](31,1)
\class[fill=red](31,1)
\class[fill=red, circlen = 3](35,1)
\class[fill=red, circlen = 3](35,1)
\class[fill=red](39,1)
\class[fill=red](39,1)
\class[fill=red](43,1)
\class[fill=red](43,1)

\foreach \i in {3,4} {
	\class[fill=red](10*\i-1, 2*\i+1)
	}

\foreach \i in {3} {
	\structline[black](10*\i-1, 2*\i+1)(10*(\i+1)-1,2*(\i+1)+1)
	}

\foreach \i in {0,1,2,3} {
\class[fill=red](26+72*\i,2)
\class[fill=red](36+72*\i,4)
	}

\structline[black](26,2)(29,7)
\structline[black](26,2)(36,4)
\structline[black](36,4)(39,9)

\foreach \i in {1,2,3,4} {
	\class[fill=blue](10*\i+72, 2*\i)
	}

\foreach \i in {1,2,3} {
	\structline[black](10*\i+72, 2*\i)(10*(\i+1)+72,2*(\i+1))
	}

\structline[black](10+72,2)(13+72,3)
\structline[black](3+72,1)(13+72,3)
\structline[black](27+72,1)(30+72,6)
\structline[black](27+72,1)(37+72,3)
\structline[black](37+72,3)(40+72,8)

\class[fill=red](2+72,2)

\class[fill=red](3+72,1)
\class[fill=red](3+72,1)
\class[fill=red](3+72,1)
\class[fill=red](7+72,1)
\class[fill=red](7+72,1)
\class[fill=red](7+72,1)
\class[fill=red](7+72,1)
\class[fill=red, circlen = 2](11+72,1)
\class[fill=red, circlen = 2](11+72,1)
\class[fill=red, circlen = 2](11+72,1)
\class[fill=red, circlen = 2](11+72,1)
\class[fill=red](15+72,1)
\class[fill=red](15+72,1)
\class[fill=red](15+72,1)
\class[fill=red](15+72,1)
\class[fill=red](19+72,1)
\class[fill=red](19+72,1)
\class[fill=red](19+72,1)
\class[fill=red](19+72,1)
\class[fill=red, circlen = 2](23+72,1)
\class[fill=red, circlen = 2](23+72,1)
\class[fill=red, circlen = 2](23+72,1)
\class[fill=red, circlen = 2](23+72,1)
\class[fill=red, circlen = 2](23+72,1)
\class[fill=red](27+72,1)
\class[fill=red](27+72,1)
\class[fill=red](27+72,1)
\class[fill=red](27+72,1)
\class[fill=red](31+72,1)
\class[fill=red](31+72,1)
\class[fill=red](31+72,1)
\class[fill=red](31+72,1)
\class[fill=red](31+72,1)
\class[fill=red, circlen = 3](35+72,1)
\class[fill=red, circlen = 3](35+72,1)
\class[fill=red, circlen = 3](35+72,1)
\class[fill=red, circlen = 3](35+72,1)
\class[fill=red, circlen = 3](35+72,1)
\class[fill=red](39+72,1)
\class[fill=red](39+72,1)
\class[fill=red](39+72,1)
\class[fill=red](39+72,1)
\class[fill=red](39+72,1)
\class[fill=red](43+72,1)
\class[fill=red](43+72,1)
\class[fill=red](43+72,1)
\class[fill=red](43+72,1)
\class[fill=red](43+72,1)

\foreach \i in {1,2,3,4} {
	\class[fill=red](10*\i-1+72, 2*\i+1)
	}

\foreach \i in {1,2,3} {
	\structline[black](10*\i-1+72, 2*\i+1)(10*(\i+1)-1+72,2*(\i+1)+1)
	}

\class[fill=red](71,1)
\structline[black](71,1)(9+72,3)
\structline[black](71,1)(74,2)

\structline[black](26+72,2)(29+72,7)
\structline[black](26+72,2)(36+72,4)
\structline[black](36+72,4)(39+72,9)

\class[fill=blue](23,5)
\structline[black](23,1,2)(23,5)
\structline[black](13,3)(23,5)
\structline[black](20,4)(23,5)

\class[red, circle](2,2)\class[red, circle](9,3)\class[red, circle](12,4)\class[red, circle](19,5)
\class[red, circle](-1,1)
\structline[black, dashed](-1,1)(2,2)
\structline[black, dashed](-1,1)(9,3)
\structline[black, dashed](2,2)(12,4)
\structline[black, dashed](9,3)(19,5)
\structline[black, dashed](9,3)(12,4)
\structline[black, dashed](19,5)(29,7)

\end{sseqdata}

\begin{sseqdata}[name = homotopyoftmfpsiat3zoomedout, Adams grading, tick step = 24,
major tick step = 24,
minor tick step = 4]    

\foreach \i in {1,2,3,4} {
	\class[fill=blue](10*\i, 2*\i)
	}

\class[fill=blue, rectangle](0,0)

\foreach \i in {0,1,2,3} {
	\structline[fill=black](10*\i, 2*\i)(10*(\i+1),2*(\i+1))
	}

\foreach \i in {0,1,2,3} {
\class[fill=blue](3+72*\i,1)
\class[fill=blue](13+72*\i,3)
\class[fill=blue](27+72*\i,1)
\class[fill=blue](37+72*\i,3)
	}

\structline[black](0,0)(3,1)
\structline[black](10,2)(13,3)
\structline[black](3,1)(13,3)
\structline[black](27,1)(30,6)
\structline[black](27,1)(37,3)
\structline[black](37,3)(40,8)

\class[fill=red](7,1)
\class[fill=red, circlen = 2](11,1)
\class[fill=red](15,1)
\class[fill=red](19,1)
\class[fill=red, circlen = 2](23,1)
\class[fill=red](23,1)
\class[fill=red](27,1)
\class[fill=red](31,1)
\class[fill=red](31,1)
\class[fill=red, circlen = 3](35,1)
\class[fill=red, circlen = 3](35,1)
\class[fill=red](39,1)
\class[fill=red](39,1)
\class[fill=red](43,1)
\class[fill=red](43,1)

\foreach \i in {3,4} {
	\class[fill=red](10*\i-1, 2*\i+1)
	}

\foreach \i in {3} {
	\structline[black](10*\i-1, 2*\i+1)(10*(\i+1)-1,2*(\i+1)+1)
	}

\foreach \i in {0,1,2,3} {
\class[fill=red](26+72*\i,2)
\class[fill=red](36+72*\i,4)
	}

\structline[black](26,2)(29,7)
\structline[black](26,2)(36,4)
\structline[black](36,4)(39,9)

\foreach \i in {1,2,3,4} {
	\class[fill=blue](10*\i+72, 2*\i)
	}

\foreach \i in {1,2,3} {
	\structline[black](10*\i+72, 2*\i)(10*(\i+1)+72,2*(\i+1))
	}

\structline[black](10+72,2)(13+72,3)
\structline[black](3+72,1)(13+72,3)
\structline[black](27+72,1)(30+72,6)
\structline[black](27+72,1)(37+72,3)
\structline[black](37+72,3)(40+72,8)

\class[fill=red](2+72,2)

\class[fill=red](3+72,1)
\class[fill=red](3+72,1)
\class[fill=red](3+72,1)
\class[fill=red](7+72,1)
\class[fill=red](7+72,1)
\class[fill=red](7+72,1)
\class[fill=red](7+72,1)
\class[fill=red, circlen = 2](11+72,1)
\class[fill=red, circlen = 2](11+72,1)
\class[fill=red, circlen = 2](11+72,1)
\class[fill=red, circlen = 2](11+72,1)
\class[fill=red](15+72,1)
\class[fill=red](15+72,1)
\class[fill=red](15+72,1)
\class[fill=red](15+72,1)
\class[fill=red](19+72,1)
\class[fill=red](19+72,1)
\class[fill=red](19+72,1)
\class[fill=red](19+72,1)
\class[fill=red, circlen = 2](23+72,1)
\class[fill=red, circlen = 2](23+72,1)
\class[fill=red, circlen = 2](23+72,1)
\class[fill=red, circlen = 2](23+72,1)
\class[fill=red, circlen = 2](23+72,1)
\class[fill=red](27+72,1)
\class[fill=red](27+72,1)
\class[fill=red](27+72,1)
\class[fill=red](27+72,1)
\class[fill=red](31+72,1)
\class[fill=red](31+72,1)
\class[fill=red](31+72,1)
\class[fill=red](31+72,1)
\class[fill=red](31+72,1)
\class[fill=red, circlen = 3](35+72,1)
\class[fill=red, circlen = 3](35+72,1)
\class[fill=red, circlen = 3](35+72,1)
\class[fill=red, circlen = 3](35+72,1)
\class[fill=red, circlen = 3](35+72,1)
\class[fill=red](39+72,1)
\class[fill=red](39+72,1)
\class[fill=red](39+72,1)
\class[fill=red](39+72,1)
\class[fill=red](39+72,1)
\class[fill=red](43+72,1)
\class[fill=red](43+72,1)
\class[fill=red](43+72,1)
\class[fill=red](43+72,1)
\class[fill=red](43+72,1)

\foreach \i in {1,2,3,4} {
	\class[fill=red](10*\i-1+72, 2*\i+1)
	}

\foreach \i in {1,2,3} {
	\structline[black](10*\i-1+72, 2*\i+1)(10*(\i+1)-1+72,2*(\i+1)+1)
	}

\class[fill=red](71,1)
\structline[black](71,1)(9+72,3)
\structline[black](71,1)(74,2)

\structline[black](26+72,2)(29+72,7)
\structline[black](26+72,2)(36+72,4)
\structline[black](36+72,4)(39+72,9)

\class[red, circle](2,2)\class[red, circle](9,3)\class[red, circle](12,4)\class[red, circle](19,5)
\class[red, circle](-1,1)
\structline[black, dashed](-1,1)(2,2)
\structline[black, dashed](-1,1)(9,3)
\structline[black, dashed](2,2)(12,4)
\structline[black, dashed](9,3)(19,5)
\structline[black, dashed](9,3)(12,4)
\structline[black, dashed](19,5)(29,7)

\end{sseqdata}

\begin{sseqdata}[name = hurewicztmfpsiat3, Adams grading, tick step = 4,
major tick step = 4,
minor tick step = 1]    

\class[fill=orange, rectangle](0,0)

\class[fill=orange](7,1)
\class[fill=orange, circlen = 2](11,1)
\class[fill=orange](15,1)
\class[fill=orange](19,1)
\class[fill=orange, circlen = 2](23,1)
\class[fill=black](27,1)
\class[fill=orange](31,1)
\class[fill=orange, circlen = 3](35,1)
\class[fill=orange](39,1)
\class[fill=orange](43,1)

\foreach \i in {1,2,3,4} {
	\class[fill=orange](10*\i, 2*\i)
	}
 
\foreach \i in {0,1,2,3} {
	\structline[black](10*\i, 2*\i)(10*(\i+1),2*(\i+1))
	}

\class[fill=orange](3,1)
\class[fill=orange](13,3)
\class[fill=orange](27,1)
\class[fill=orange](37,3)

\foreach \i in {1,2,3} {
\class[fill=black](3+72*\i,1)
\class[fill=black](27+72*\i,1)
	}

 \foreach \i in {1,2,3} {
\class[fill=orange](13+72*\i,3)
\class[fill=orange](37+72*\i,3)
	}

\structline[black](0,0)(3,1)
\structline[black](10,2)(13,3)
\structline[black](3,1)(13,3)
\structline[black](27,1)(30,6)
\structline[black](27,1)(37,3)
\structline[black](37,3)(40,8)

\class[fill=orange](29,7)
\class[fill=orange](39,9)

\foreach \i in {3} {
	\structline[black](10*\i-1, 2*\i+1)(10*(\i+1)-1,2*(\i+1)+1)
	}

\class[fill=orange](26,2)
\class[fill=orange](36,4)

\structline[black](26,2)(29,7)
\structline[black](26,2)(36,4)
\structline[black](36,4)(39,9)

\foreach \i in {1,2,3,4} {
	\class[fill=orange](10*\i+72, 2*\i)
	}

\foreach \i in {1,2,3} {
	\structline[black](10*\i+72, 2*\i)(10*(\i+1)+72,2*(\i+1))
	}

\structline[black](10+72,2)(13+72,3)
\structline[black](3+72,1)(13+72,3)
\structline[black](27+72,1)(30+72,6)
\structline[black](27+72,1)(37+72,3)
\structline[black](37+72,3)(40+72,8)

\class[fill=black](3+72,1)

\class[fill=orange](12+72,4)
\class[fill=orange](2+72,2)
\class[fill=black](26+72,2)
\class[fill=orange](36+72,4)

\foreach \i in {1,2,3,4} {
	\class[fill=orange](10*\i-1+72, 2*\i+1)
	}

\foreach \i in {1,2,3} {
	\structline[black](10*\i-1+72, 2*\i+1)(10*(\i+1)-1+72,2*(\i+1)+1)
	}

\class[red](71,1)
\structline[black](71,1)(9+72,3)
\structline[black](71,1)(74,2)

\structline[black](9+72,3)(12+72,4)
\structline[black](2+72,2)(12+72,4)
\structline[black](26+72,2)(29+72,7)
\structline[black](26+72,2)(36+72,4)
\structline[black](36+72,4)(39+72,9)

\class[fill=orange](3+72,1)
\class[fill=orange](7+72,1)
\class[fill=orange, circlen = 2](11+72,1)
\class[fill=orange](15+72,1)
\class[fill=orange](19+72,1)
\class[fill=orange, circlen = 2](23+72,1)

\class[fill=black](27+72,1)
\class[fill=orange](27+72,1)
\class[fill=orange](31+72,1)
\class[fill=orange, circlen = 3](35+72,1)
\class[fill=orange](39+72,1)
\class[fill=orange](43+72,1)

\foreach \i in {1,2,3,4} {
	\class[fill=orange](10*\i+72+72, 2*\i)
	}

\foreach \i in {1,2,3} {
	\structline[black](10*\i+72+72, 2*\i)(10*(\i+1)+72+72,2*(\i+1))
	}

\structline[black](10+72+72,2)(13+72+72,3)
\structline[black](3+72+72,1)(13+72+72,3)
\structline[black](27+72+72,1)(30+72+72,6)
\structline[black](27+72+72,1)(37+72+72,3)
\structline[black](37+72+72,3)(40+72+72,8)

\class[fill=black](3+72+72,1)

\class[fill=orange](12+72+72,4)
\class[fill=black](2+72+72,2)
\class[fill=orange](26+72+72,2)
\class[fill=orange](36+72+72,4)

\class[diamond, fill=green](9+72+72,3)
\class[fill=orange](19+72+72,5)

\foreach \i in {3,4} {
	\class[fill=orange](10*\i-1+72+72, 2*\i+1)
	}

\foreach \i in {1,2,3} {
	\structline[black](10*\i-1+72+72, 2*\i+1)(10*(\i+1)-1+72+72,2*(\i+1)+1)
	}

\class[red](71+72,1)
\structline[black](71+72,1)(9+72+72,3)
\structline[black](71+72,1)(74+72,2)

\structline[black](9+72+72,3)(12+72+72,4)
\structline[black](2+72+72,2)(12+72+72,4)
\structline[black](26+72+72,2)(29+72+72,7)
\structline[black](26+72+72,2)(36+72+72,4)
\structline[black](36+72+72,4)(39+72+72,9)

\class[fill=orange](3+72+72,1)
\class[fill=orange](7+72+72,1)
\class[fill=orange, circlen = 2](11+72+72,1)
\class[fill=orange](15+72+72,1)
\class[fill=orange](19+72+72,1)
\class[fill=orange, circlen = 2](23+72+72,1)

\class[black](27+72+72,1)
\class[fill=orange](27+72+72,1)
\class[fill=orange](31+72+72,1)
\class[fill=orange, circlen = 3](35+72+72,1)
\class[fill=orange](39+72+72,1)
\class[fill=orange](43+72+72,1)

\classoptions["{\al_1}" {above left, pin=black}](3,1)
\classoptions["{\al_2}" {below, pin=black}](7,1)
\classoptions["{\al_{3/2}}" {below, pin=black}](11,1)
\classoptions["{\al_{6/2}}" {below, pin=black}](23,1)
\classoptions["{\al_{9/3}}" {below, pin=black}](35,1)

\classoptions["{\al_1\be_{3/3}}" {right, pin=black}](37,3)

\classoptions["{\be_1}" {right, pin=black}](10,2)
\classoptions["{\be_2}" {above left, pin=black}](26,2)

\classoptions["{\ga_1}" {above left, pin=black}](29,7)

\classoptions["{\be_{6/3}}" {right, pin=black}](82,2)
\classoptions["{\al_1\be_7}" {right, pin=black}](109,3)

\classoptions["{\be_5}" {above left, pin=black}](74,2)
\classoptions["{x_{81}=\langle\be_5,\al_1,\al_1\rangle}" {above left, pin=black}](81,3)
\classoptions["{\be_2\be_{6/3}}" {above left, pin=black}](108,4)

\classoptions["{\be_5\be_{6/3}}" {above, pin=black}](12+144,4)
\classoptions["{x_{81}\be_{6/3}}" {above, pin=black}](163,5)
\classoptions["{\be_{10}}" {right, pin=black}](10+144,2)
\classoptions["{\be_{11}}" {left, pin=black}](26+144,2)
\classoptions["{\al_1\be_{12/3}}" {below, pin=black}](37+144,3)
\classoptions["{\be_{6/3}^2=\be_1\be_{10}}" {below, pin=black}](20+144,4)
\classoptions["{\langle \beta_{9/4}, 3, \be_2\rangle}" {above left, pin=black}](9+144,3)

\class[fill=black](23,1)
\class[fill=orange](23,5)
\structline[black](13,3)(23,5)
\structline[black](20,4)(23,5)
\structline[black](23,1,2)(23,5)

\class[fill=black](23+72,1)
\class[fill=orange](23+72,5)
\structline[black](13+72,3)(23+72,5)
\structline[black](20+72,4)(23+72,5)
\structline[black](23+72,1,2)(23+72,5)

\class[fill=black](23+144,1)
\class[fill=orange](23+144,5)
\structline[black](13+144,3)(23+144,5)
\structline[black](20+144,4)(23+144,5)
\structline[black](23+144,1,2)(23+144,5)

\end{sseqdata}

\begin{sseqdata}[name = homotopyoftmfmod3, Adams grading, tick step = 4,
major tick step = 4,
minor tick step = 1]    

\class[fill=blue](0,0)
\class[fill=blue](8,0)
\class[fill=blue](12,0)
\class[fill=blue](16,0)
\class[fill=blue](20,0)
\class[fill=blue](24,0)
\class[fill=blue](24,0)
\class[fill=blue](28,0)
\class[fill=blue](32,0)
\class[fill=blue](32,0)
\class[fill=blue](36,0)
\class[fill=blue](36,0)
\class[fill=blue](40,0)
\class[fill=blue](40,0)

\class[fill=blue](3,1)
\class[fill=blue](13,3)

\class[fill=blue](10,2)
\class[fill=blue](20,4)
\class[fill=blue](30,6)
\class[fill=blue](40,8)

\class[fill=blue](27,1)
\class[fill=blue](37,3)

\class[fill=red](4,0)
\class[fill=red](14,2)

\class[fill=red](11,1)
\class[fill=red](21,3)
\class[fill=red](31,5)
\class[fill=red](41,7)

\class[fill=red](28,0)
\class[fill=red](38,2)

\structline[black](0,0)(3,1)
\structline[black](0,0)(10,2)
\structline[black](3,1)(13,3)
\structline[black](10,2)(13,3)

\structline[black](10,2)(20,4)
\structline[black](20,4)(30,6)
\structline[black](30,6)(40,8)

\structline[black](27,1)(30,6)
\structline[black](37,3)(40,8)
\structline[black](27,1)(37,3)

\structline[black](4,0)(14,2)
\structline[black](11,1)(14,2)
\structline[black](11,1)(21,3)
\structline[black](21,3)(31,5)
\structline[black](31,5)(41,7)

\structline[orange](0,0)(4,0)
\structline[orange](4,0)(8,0)
\structline[orange](8,0)(12,0)
\structline[orange](12,0)(16,0)
\structline[orange](16,0)(20,0)
\structline[orange](20,0)(24,0,1)

\structline[orange](32,0)(36,0)
\structline[orange](36,0)(40,0)

\structline[orange](10,2)(14,2)

\classoptions["{\overline{\al}}" above left](3,1)
\classoptions["{\overline{\be}}" above left](10,2)
\classoptions["{v_1}" right](4,0)

\end{sseqdata}

\begin{sseqdata}[name = purejtwohomotopy, Adams grading, tick step = 4,
major tick step = 4,
minor tick step = 1]    

\foreach \i in {1,2,3,4} {
	\class[fill=blue](10*\i, 2*\i)
	}

\class[fill=blue, rectangle](0,0)

\foreach \i in {0,1,2,3} {
	\structline[fill=black](10*\i, 2*\i)(10*(\i+1),2*(\i+1))
	}

\foreach \i in {0,1,2,3} {
\class[fill=blue](3+72*\i,1)
\class[fill=blue](13+72*\i,3)
\class[fill=blue](27+72*\i,1)
\class[fill=blue](37+72*\i,3)
	}

\structline[black](0,0)(3,1)
\structline[black](10,2)(13,3)
\structline[black](3,1)(13,3)
\structline[black](27,1)(30,6)
\structline[black](27,1)(37,3)
\structline[black](37,3)(40,8)

\class[fill=red](7,1)
\class[fill=red, circlen = 2](11,1)
\class[fill=red](15,1)
\class[fill=red](19,1)
\class[fill=red, circlen = 2](23,1)
\class[fill=red, circlen = 2](23,1)
\class[fill=red](27,1)
\class[fill=red](31,1)
\class[fill=red](31,1)
\class[fill=red, circlen = 3](35,1)
\class[fill=red, circlen = 3](35,1)
\class[fill=red](39,1)
\class[fill=red](39,1)
\class[fill=red](43,1)
\class[fill=red](43,1)

\foreach \i in {1,2,3,4} {
	\class[fill=red](10*\i-1, 2*\i+1)
	}

\foreach \i in {1,2,3} {
	\structline[black](10*\i-1, 2*\i+1)(10*(\i+1)-1,2*(\i+1)+1)
	}

\foreach \i in {0,1,2,3} {
\class[fill=red](12+72*\i,4)
\class[fill=red](26+72*\i,2)
\class[fill=red](36+72*\i,4)
	}

\structline[black](9,3)(12,4)
\structline[black](26,2)(29,7)
\structline[black](26,2)(36,4)
\structline[black](36,4)(39,9)
\end{sseqdata}

\begin{sseqdata}[name = purejLtauhomotopy, Adams grading, tick step = 4,
major tick step = 4,
minor tick step = 1]    

\foreach \i in {1,2,3,4} {
	\class[fill=blue](10*\i, 2*\i)
	}

\class[fill=blue, rectangle](0,0)

\foreach \i in {0,1,2,3} {
	\structline[fill=black](10*\i, 2*\i)(10*(\i+1),2*(\i+1))
	}

\foreach \i in {0,1,2,3} {
\class[fill=blue](3+72*\i,1)
\class[fill=blue](13+72*\i,3)
\class[fill=blue](27+72*\i,1)
\class[fill=blue](37+72*\i,3)
	}

\structline[black](0,0)(3,1)
\structline[black](10,2)(13,3)
\structline[black](3,1)(13,3)
\structline[black](27,1)(30,6)
\structline[black](27,1)(37,3)
\structline[black](37,3)(40,8)

\class[fill=red](7,1)
\class[fill=red, circlen = 2](11,1)
\class[fill=red](15,1)
\class[fill=red](19,1)
\class[fill=red, circlen = 2](23,1)
\class[fill=red, circlen = 2](23,1)
\class[fill=red](27,1)
\class[fill=red](31,1)
\class[fill=red](31,1)
\class[fill=red, circlen = 3](35,1)
\class[fill=red, circlen = 3](35,1)
\class[fill=red](39,1)
\class[fill=red](39,1)
\class[fill=red](43,1)
\class[fill=red](43,1)

\foreach \i in {3,4} {
	\class[fill=red](10*\i-1, 2*\i+1)
	}

\foreach \i in {3} {
	\structline[black](10*\i-1, 2*\i+1)(10*(\i+1)-1,2*(\i+1)+1)
	}

\foreach \i in {0,1,2,3} {
\class[fill=red](26+72*\i,2)
\class[fill=red](36+72*\i,4)
	}

\structline[black](26,2)(29,7)
\structline[black](26,2)(36,4)
\structline[black](36,4)(39,9)

\class(4,-1)\class(34,9)
\structline[green](4,-1)(34,9)
\end{sseqdata}

\begin{document}
\title{
Nonvanishing of products in \texorpdfstring{$v_2$}{v2}-periodic families\\
at the prime 3
}
\author{Christian Carrick\footnote{\href{mailto:carrick@math.uni-bonn.de}{\texttt{carrick@math.uni-bonn.de}}}\, and Jack Morgan Davies\footnote{\href{mailto:davies@uni-wuppertal.de}{\texttt{davies@uni-wuppertal.de}}}}

\date{\today}

\maketitle

\begin{abstract}
    Many products amongst $v_2$-periodic families in the stable homotopy groups of spheres are shown not to vanish. This is done by carefully studying the action of Adams operations on topological modular forms. 
    A crucial ingredient is Pstragowski's category of synthetic spectra which affords us the necessary freedom to work with (modified) Adams--Novikov spectral sequences.
\end{abstract}

\setcounter{tocdepth}{2}
\tableofcontents

\section{Introduction}
This paper studies the homotopy groups of the \emph{height $2$ connective image-of-$J$ spectrum} $\j^2$ at the prime $3$, as defined by the second-named author. We show that $\j^2$ is an effective connective model for a particular slice of the $\K(2)$-local sphere, balancing the richness and complications of GHMR resolutions \cite{ktwolocalsphere} and Behrens' $Q(2)$ spectrum \cite{ktwospheremark} with the sparseness and simplicity of the height 1 image-of-$J$ spectrum of Adams and Mahowald.\\

We compute the image of $\pi_\ast \Sph$ inside $\pi_\ast \j^2$, also known as the \emph{Hurewicz image} of $\j^2$; see \Cref{maintheorem}. As a consequence of the multiplicative structure of $\j^2$ and its Hurewicz image, we deduce the nonvanishing of a large class of products of $v_2$-periodic families in $\pi_*\Sph_{(3)}$; see \Cref{sphereconsequences}. Notably, these $v_2$-periodic families lie in arbitrarily high stem, so beyond the finite range of complete computations of $\pi_\ast \Sph_{(3)}$, and these products also lie in Adams--Novikov filtration $\geq 4$, beyond the $E_2$-page computations of Miller--Ravenel--Wilson \cite{MRW77}.\\

The homotopy groups of $\j^2$ are easy to calculate, but detection statements such as \Cref{maintheorem} require much more input. As computing the Adams or Adams--Novikov spectral sequence for $\j^2$ seems like a difficult task, we instead define a \emph{synthetic spectrum} (or a $\mathbf{C}$-motivic spectrum or filtered spectrum, if the reader prefers) $\j^2_\BP$ using the \emph{vertical $t$-structures} of \cite[\textsection2]{syntheticj} providing us with an easy-to-compute \emph{modified Adams--Novikov spectral sequence} for $\j^2$.

\subsection*{Motivation}
Chromatic homotopy organizes the stable homotopy groups of spheres $\pi_*\Sph$ into periodic families. This began with Adams' work on the image of the $J$-homomorphism \cite{adamsjofx}
\[J\colon \pi_{*}O\to\pi_{*}\Sph\]
which unveils some $8$-fold periodic phenomenæ in $\pi_\ast \Sph$ courtesy of Bott periodicity. Adams gave another interpretation of this periodicity by constructing self-maps $v_1\colon \Sph/p[d] \to \Sph/p$ of the mod $p$ Moore spectrum for each prime $p$, where $d=8$ if $p=2$ and $d=2p-2$ otherwise. Using the maps in the cofibre sequence $\Sph\xrightarrow{p}\Sph\xrightarrow{q_0}\Sph/p\xrightarrow{\partial_0}\Sph[1]$, one defines the class $\alpha_n$ as the composite
\begin{equation}\label{alphacomposite} \al_n \colon \Sph[dn] \xrightarrow{q_0} \Sph/p[dn] \xrightarrow{v_1} \Sph/p[d(n-1)] \xrightarrow{v_1} \cdots \xrightarrow{v_1} \Sph/p \xrightarrow{\partial_0} \Sph[1]
\end{equation}
which can be used to describe the image of $J$. This perspective allows one to view the periodicity in the image of $J$ as the $v_1$-periodicity of the family of elements $\{\al_n\}$.\\

This approach to constructing periodic families in $\pi_\ast \Sph_{(p)}$ was then generalised by Smith \cite{smithbeta}, who defined a $v_2$-periodic family known as the $\beta$-family at primes $p\ge 5$. In their landmark paper, Miller--Ravenel--Wilson \cite{MRW77} constructed the chromatic spectral sequence which shows that, in a suitably attenuated sense, \emph{everything} in $\pi_*\Sph_{(p)}$ comes from such periodic families. Moreover, these families are stratified by \emph{height}, an invariant in the theory of formal groups, whereby one speaks of $v_h$-periodic families with corresponding height $h$. Adams' $\alpha$-family corresponds to height $1$, Smith's $\beta$-family to height $2$, and Miller--Ravenel--Wilson define further ``Greek letter families'' at higher heights.\\

Definitions of Greek letter families such as (\ref{alphacomposite}) are easy to make with these $v_n$-maps at hand, but it is another matter to show that these elements are nonzero. The \emph{Adams--Novikov spectral sequence} (ANSS) gives one method of verifying this. Roughly speaking, the Adams--Novikov spectral sequence arises from a filtration of the sphere $\Sph_{(p)}$ given by placing $v_n$-periodic families in filtration $n$. Using their chromatic spectral sequence, Miller--Ravenel--Wilson obtained complete information about the $E_2$-page of the ANSS in filtrations $\le 2$, from which the nontriviality of the $\al$-family follows easily, for example. Similarly, the $v_2$-periodic $\be$-families of Smith (at primes $\geq 5$) and Behrens--Pemmaraju \cite{marksatya} (at $p=3$) are detected in AN-filtration 2. A quick look at the $E_2$-page computations of Miller--Ravenel--Wilson shows these families are nonzero and in too low of a filtration to be hit by a differential in this spectral sequence; hence these $\be$-families are also nontrivial.\\

This approach to verifying the nontriviality of a family becomes much more subtle when both the stem and Adams--Novikov filtration of an element are very large. For example, using the argument above the elements $\beta_{10}$, $\beta_{19}$, and $\beta_{23}$ are nontrivial in $\pi_{*}\Sph_{(3)}$, but their product $\beta_{10}\beta_{19}\beta_{23}$ lives in the $814$ stem and Adams--Novikov filtration $\ge 6$. This is beyond the range of the computations of Miller--Ravenel--Wilson or a direct stem-by-stem approach via the Adams spectral sequence (ASS) or the ANSS. In these filtrations, there is also lots of noise generated by $v_h$-periodic families for $h\geq 3$, of which hardly any computations are known. For these simple reasons, the nontriviality of many such simple products can be mysterious.\\

The ideal situation to access such products would be to have a connective ring spectrum $E$ that detects the above families and whose homotopy groups are extremely sparse, consisting mostly of products and brackets formed from this family. In this case, it is easy to check if one of these products is nonzero in $\pi_\ast E$, and it follows immediately that it must also be nonzero in $\pi_\ast \Sph$. The $\K(2)$-local sphere $\Sph_{\K(2)}$ definitely captures all of the desired family, but the multiplicative structure of its homotopy groups is a mystery. This stems from the fact that these groups are incredibly difficult to compute. Indeed, the GHMR resolutions of \cite{ktwolocalsphere} are an attempt to reinterpret the work of Shimomura et al. \cite{l2atprime3}, \cite{ltwomooreatthree}, \cite{shimomurawang}. Even at primes $p\geq 5$, when the ANSS for $\Sph_{\K(2)}$ collapses, this computation is far from easy; see \cite{shimomurayabe}. Behrens' spectrum $Q(2)$ of \cite{ktwospheremark} captures exactly half of the $\K(2)$-local sphere in a certain sense, which substantially simplifies computations. However, there is a drastically simpler spectrum that detects the same amount of products of $v_2$-periodic families as $Q(2)$.\\

The \emph{height 2 connective image-of-$J$ spectrum} $\j^2$, introduced by the second-named author in \cite{adamsontmf}, gives us the ideal situation described above. Indeed, the homotopy groups of $\j^2$ are easily computed from the homotopy groups of $\tmf_{3}$ and the action of the Adams operation $\psi^2$ on these homotopy groups. From the perspective of chromatic homotopy theory, the $\K(2)$-localisation of $\tmf$ is the $G_{24}$-fixed points of a height 2 Lubin--Tate theory, where $G_{24}$ is the binary tetrahedral group and the maximal finite subgroup of the Morava stabiliser group $\mathbf{G}_2$. As shown in \cite[Th.6.5]{tmfthree}, these fixed points pick out the $\be_{1+9t}$- and $\be_{6+9t/3}$-families. The $\K(2)$-localisation of $\j^2$ is the further fixed points with respect to the subgroup of $\mathbf{G}_2$ generated by the Adams operation $\psi^2$. From this, one might suspect that $\j^2$ should detect more $v_2$-periodic families than $\tmf$. The real challenge undertaken in this article is confirming this suspicion by computing the image of $\pi_\ast \Sph$ inside $\pi_\ast \j^2$.

\subsection*{Main results}
The construction here of $\j^2$ at the prime $3$ is closely related to the \emph{height $1$ connective image-of-$J$ spectrum} $\j^1$ at the prime 2. One can define $\j^1$ by first setting $\ko^\psi = \ko_2^{h\Z}$ to be the homotopy fixed points of $\ko_2$ with $\Z$-action generated by the Adams operation $\psi^3$, ie, the fibre of $\psi^3-1$. Then one defines $\j^1$ as the pullback of the span 
\[\tau_{\geq 0}\ko^\psi \to \tau_{\leq 2}\tau_{\geq 0}\ko^\psi \gets \tau_{\leq 2} \Sph_2.\]
This extra round of covers and truncations cuts down on the noise in the homotopy groups of $\ko^\psi$. The result is that the unit map $\Sph \to \j^1$ induces an isomorphism on homotopy groups in low degrees and is surjective in all degrees; see \cite[Th.A]{syntheticj}. Similarly, we define $\j^2$ as the pullback of the span
\[\tau_{\geq 0}\tmf^\psi \to \tau_{\leq 22}\tau_{\geq 0}\tmf^\psi \gets \tau_{\leq 22} \Sph_3\]
where $\tmf^\psi = \tmf_3^{h\Z}$ where $\Z$ acts on $\tmf_3$ through the Adams operation $\psi^2$ of \cite[Th.B]{adamsontmf}. From this construction, we obtain a commutative diagram of $\E_\infty$-rings
\[\begin{tikzcd}
    {\j^2}\ar[r]\ar[d]  &   {\tmf_3}\ar[d]  \\
    {\j^1}\ar[r]        &   {\ko_3,}
\end{tikzcd}\]
where by abuse of notation, we now write $\j^1$ for the \emph{height 1 connective image-of-$J$ spectrum} at the prime $3$, related to the fibre of $\psi^2-1$ acting on $\ko_3$. It follows that all elements in $\pi_\ast \Sph$ detected by $\j^1$ or $\tmf_3$ are also detected by $\j^2$. The second-named author asked in \cite[Q.4.12]{adamsontmf}, if this is a complete description of the Hurewicz image of $\j^2$. One salient feature of $\j^2$ is its ability to detect many more elements in $\pi_\ast \Sph$.

\begin{theoremalph}\label{maintheorem}
The Hurewicz image of $\j^2$ is given by the disjoint union of the families\footnote{This $v_2^9$-periodic family $x_{81}^{(t)}$ is defined in \Cref{df:newperiodicfamily} and is generated by the element $x_{81}^{(0)} = x_{81} = \langle \al_1, \al_1, \beta_5\rangle$.}
\begin{equation}\label{hurewiczofj}
    \al_{a/\nu_3(a)+1}
\end{equation}
\begin{equation}\label{hurewiczoftmf}
    \be_1^k\be_{1+9t}, \quad 
    \al_1\be_{1+9t}, \quad 
    [\al_1\be_{3+9t/3}], \quad
    \be_1^k\be_{6+9t/3}, \quad 
    \al_1\be_{6+9t/3}, \quad 
    [\al_1\be_{7+9t}]
\end{equation}
\begin{equation}\label{hurewiczofjtwoa} 
    \al_1\be_1\be_{1+9t}, \quad
    \be_1[\al_1\be_{3+9t/3}], \quad 
    \al_1\be_1\be_{6+9t/3},  \quad 
    \be_1[\al_1\be_{7+9t}],
\end{equation}
\begin{equation}\label{hurewiczofjtwob}
\al_1^i\be_1^j \be_{2+9t}, \quad
    \be_{6/3} \be_{2+9t}, \quad
    \be_1^j \be_{5+9t}, \quad
    \be_{6/3} \be_{5+9t}, \quad 
    \be_1^k x_{81}^{(t)}, \quad 
    \be_{6/3}x_{81}^{(t)}
\end{equation}
for $1\leq a$, $0\leq i,j\leq 1$, $0\leq k\leq 3$, and $0\leq t$, with the exception of the classes in $\pi_{9+144s} \j^2$ for $1\leq s$, which are ambiguous (\Cref{ambiguousclasses,153question}). This Hurewicz image is depicted by the orange and yellow classes in \Cref{hurewiczofjpartone,hurewiczofjparttwo} with the exception in green.
\end{theoremalph}


The first family of elements (\ref{hurewiczofj}) comes from $\j^1$ and the second (\ref{hurewiczoftmf}) from $\tmf_3$; the third (\ref{hurewiczofjtwoa}) and fourth (\ref{hurewiczofjtwob}) families are new to $\j^2$. This answers \cite[Q.4.12]{adamsontmf} in the positive: there are classes detected by $\j^2$ which are not a simple combination of those from $\j^1$ and $\tmf_3$.\\

The description of the Hurewicz image of $\j^2$ above is \emph{minimal}, in other words, we have described elements in $\pi_\ast\Sph$ which map to generators of the Hurewicz image of $\j^2$. There are however, more elements in $\pi_\ast \Sph$ which are detected by $\j^2$ given by various products. For example, inside $\j^2$, we know that the classes $\be_1\be_{11}$ and $\be_{10}\be_2$ are detected by the same nonzero class in $\pi_\ast \j^2$ (up to a unit), but only $\be_1\be_{11}$ is described in (\ref{hurewiczofjtwoa}). A more extensive list of these products appears in \Cref{listofproducts}. Employing this larger collection of elements in $\pi_\ast \Sph$ which are detected in $\j^2$ as well as some Toda bracket arguments, we obtain the following nonvanishing statements in the stable homotopy groups of spheres.

\begin{theoremalph}\label{sphereconsequences}
    Let $0\leq s,t,w$ and $x_i \in \{\be_{1+9s}, \be_{6+9s/3}\}$ be a collection of elements. Then the following products are nonzero in $\pi_\ast\Sph$:
    \begin{equation}\label{alternativeproductstmf}
        \prod_{i=1}^4 x_i, \qquad \al_1(\al_1\be_{3+9s/3}), \qquad \al_1(\al_1\be_{7+9s}), \qquad 
        \al_1\prod_{i=1}^2 x_i, \qquad
        \al_1 x_i \be_{2+9t}, \qquad 
        x_i\be_{5+9t}\end{equation}
    \begin{equation}\label{alternativeproductsjtwoa}
        \al_1 x_{81}^{(s)}, \quad x_{81}^{(s)} \prod_{i=1}^3 x_i, \quad x_i [\al_1\be_{3+9s/3}], \quad x_i [\al_1\be_{7+9s}]
    \end{equation}
    \begin{equation}\label{jnondetectioncor} \be_{1+9s}\be_{1+9t}\be_{5+9w}, \qquad \be_{6+9s/3}\be_{6+9t/3}\be_{5+9w}\end{equation}
\end{theoremalph}

The nonvanishing of the families (\ref{alternativeproductstmf}) and (\ref{alternativeproductsjtwoa}) all follow immediately from the fact that these products are nonzero in $\j^2$. On the other hand, the nonvanishing of the families in (\ref{jnondetectioncor}) follows from the fact that certain classes in $\pi_\ast\j^2$ do \textbf{not} lie in the Hurewicz image; see \Cref{tmftodacorollary,todacorolalry}. As a consequence of \Cref{sphereconsequences}, we also obtain nonvanishing results for elements in $\pi_\ast \Sph/3$; see \Cref{sillyaprime,moorespectrumnonvansihing}. These results simultaneously recover and generalise various results of Arita--Shimomura \cite{aritashimomuromodthreemoore} and Shimomura \cite{shimbetaoneactionatthree}; see \Cref{generalisationremark}. These are further built upon by the second-named author in \cite{davies_beta}.\\

In particular, we obtain a large number of new $v_2^9$-periodic families inside $\pi_\ast \Sph$ at the prime $3$.\\

The main tool we use to study $\j^2$ is \emph{synthetic spectra}. As far as the techniques of this article are concerned, one could alternatively use the \emph{$\mathbf{C}$-motivic modular forms} of \cite{mmf}. The $\infty$-category of synthetic spectra $\Syn$,  introduced by Pstr\k{a}gowski \cite{syntheticspectra}, can be thought of as an $\infty$-category of ANSSs. For example, there is a functor $\nu \colon \Sp \to \Syn$ such that $\nu X$ exactly encodes the ANSS of $X$. In particular, $\nu \tmf_3$ encodes the well-known ANSS for $\tmf_3$ as described in \cite{bauer} or \cite{smfcomputation}, for example. We then define a synthetic spectrum $\j^2_\BP$ that encodes a \emph{modified} ANSS for $\j^2$ and is also readily computable from the fibre of the map $\psi^2-1 \colon \nu\tmf \to \nu\tmf$ acting on the ANSS for $\tmf$. As mentioned in \cite{syntheticj}, the classical ANSS for $\j^1$ does not appear in the literature, let alone an ANSS for $\j^2$. These synthetic methods provide us with a viable inroad to detection statements regarding $\j^2$, eventually leading us to \Cref{maintheorem}.\\

This article is a sequel to \cite{syntheticj}, where we defined linear $t$-structures on $\Syn$ and used this to provide a simple proof of the classical statement that $\j^1$ detects the divided $\al$-family. Together, these two articles are the first steps in a wider project of detection statements of ring spectra using modern techniques of synthetic spectra, another instance of this being the fundamental detection statements for $\tmf$ in \cite[\textsection5]{smfcomputation}. A quick modified AN-filtration argument (similar to \Cref{filtrationone}) shows that $\j^2$ at $p\geq 5$ can only detect the divided $\alpha$-family, but at $p=2$, $\j^2$ can be used to detect many new $v_2^{32}$-periodic families in the stable stems; see \cite{v232families} as well as applications to the geometry of exotic spheres in \cite{bauer_quigley}. In \cite{davies_beta}, the second-named author has also refined some of the computations here, generalising Shimomura's work \cite{shimbetaoneactionatthree}. We are also exploring forms of Behrens' $Q(N)$ spectra of \cite{ktwospheremark}; see \Cref{qnremark}. Equivariant versions of $\J^2$ and $Q(2)$ will also appear in forthcoming work of Balderrama--Linskens with the second-named author \cite{geometricnorms}.


\subsection*{Outline}

In \Cref{homotopygroupssection}, we (re)define $\j^2$ (\Cref{definitionofjtwonew}), compute its homotopy groups (see \Cref{ss1,ss2,ss3,ss4}), and define a preferred synthetic lift $\j^2_\BP$ (\Cref{syntheticjdefinition}) as well as periodic variants $\J^2$ and $\J^2_\BP$ (\Cref{periodicdefinition}).\\

In \Cref{detectionsection}, we prove \Cref{maintheorem}. This begins with the formal arguments that the Hurewicz images of $\tmf$ and $\j^1$ are contained in that of $\j^2$ (\Cref{hurewzicoftmf,alphafamilydetection}). Next, we directly prove, using the synthetic detection techniques of \cite{syntheticj}, that $\j^2$ detects the $\beta_2$-family (\Cref{betatwofamily}) and from the multiplicative structure of $\j^2$ we also obtain the whole Hurewicz image in degree $\leq 144$ (\ref{lowdetection}). Then we use $v_2^9$-periodic phenomen{\ae} in $\j^2$ (\Cref{v29periodicity}), similar to the $v_1$-periodic Adams periodicity operator $\langle \sigma, 16, -\rangle$ on $\pi_\ast \Sph$, to produce a lower bound on the Hurewicz image of $\j^2$ (\Cref{periodicconclusions}). Finally, we give a series of ad hoc results concerning the nondetection of many elements in $\pi_\ast\j^2$ (\Cref{filtrationone,filtrationtwonondetection}), thus completing the proof of \Cref{maintheorem}.\\

In \Cref{nonvanishingsection}, we prove \Cref{sphereconsequences} by combining the detection statements of \Cref{maintheorem} together with the multiplicative structure on $\j^2$. We also mention some applications to the homotopy groups of the Moore spectrum $\mathbf S/3$ in \Cref{ssec:Moore}.


\subsection*{Notation}
From now on, all spectra will be implicitly $3$-completed. We will also assume the reader has a familiarity with the homotopy groups of $\tmf$ at the prime $3$; see \cite[\textsection5-6]{bauer}, \cite[\textsection13]{tmfbook}, or \cite[\textsection7.1]{smfcomputation} for details. We write elements in $\pi_\ast \tmf_3$ using their $E_2$-page representative and no other decoration. For example, a generator of $\pi_{27}\tmf_3\simeq \F_3$ will be written as $\al\Delta$ even though this class is \textbf{not} $\al$-divisible in homotopy groups. We write $\al$ and $\be$ with no decoration of the associated classes in $\pi_\ast \tmf$; in $\pi_\ast \j^2$ we decorate these classes as $\al_1$ and $\be_1$, respectively. We reserve the symbol $\partial$ for the boundary map $\tmf[-1] \to \tmf^\psi$ from the cofibre sequence (\ref{fibredefinition}) or its synthetic version (\ref{definitingcofibresynthic}).\\

We write $\ueq$ for equality up to a unit in $\F_3$. Given a spectrum $X$ and the cofibre sequence
\[X\xrightarrow{ 3^r\cdot} X \xrightarrow{q_{0^r}} X/3^r \xrightarrow{\partial_{0^r}} X[1],\]
we write $\overline{x}=q_0(x)$ for any $x\in \pi_d X$ and $\widetilde{x}\in \pi_d X/3$ for a choice of element such that $\partial_0(\widetilde{x})=x\in \pi_{d-1}X$; of course, this class $\widetilde{x}$ is not generally uniquely determined. For typographical reasons, we will also write $\widetilde{x} = (x)^\sim$. We implicitly use the fact that each $\pi_d X/3$ is a $\F_3$-vector space as $\Sph/3$ has a unital multiplication by \cite[Ex.2.2]{okaringstructures} and $\pi_0\Sph/3\simeq \F_3$. We write $v_1, v_2$, and their powers for both the elements in homotopy groups and the associated self-maps---it will be clear what we mean by context. We define the $\be_{i/j}$ elements, if they exist, in $\pi_{4(4i-j)-2}\Sph$ as $\partial_0\partial_{1^j}(v_2^i)$, where the boundary maps are associated with the cofibre sequences
\[\Sph \xrightarrow{3} \Sph \xrightarrow{q_0} \Sph/3\xrightarrow{\partial_0} \Sph[1]  \qquad \qquad \Sph/3\xrightarrow{v_1^j} \Sph/3 \xrightarrow{q_{1^j}} \Sph/(3,v_1^j) \xrightarrow{\partial_{1^j}} \Sph[4j+1]/3\]
and $v_2^i$ is an element in $\pi_{16i}\Sph/(3,v_1^j)$ detected by the element $v_2^i$ from the $E_2$-page of the ANSS.\\

Finally, we will freely use the language of ($\BP$-based) \emph{synthetic spectra} as defined in \cite{syntheticspectra} and using the notation of \cite{syntheticj}. In particular, we use ``stem--filtration'' grading for $\Syn=\Syn_{\BP}$, so $\pi_{k,s}$ corresponds to a $(k,s)$-location in an Adams chart. Formally, this means that we write $\Sph[a,b]:=(\nu\Sph[a+b])[-b]$ and $X[a,b] = \Sph[a,b]\otimes X$. Write $\1 = \Sph[0,0]$ for the unit. In particular, the class $\tau \in \pi_{\ast,\ast} \1$ has degree $(0,-1)$. Using Pstr\k{a}gowski's equivalence of \cite[Th.1.4]{syntheticspectra}, one could also replace all of the $3$-complete synthetic spectra (which happen to be even) with $3$-complete \emph{$\mathbf{C}$-motivic spectra}, if the reader prefers. One could also work solely in the language of modules over the ANSS for $\Sph$ in filtered spectra.


\subsection*{Acknowledgements}
Thank you to Eva Belmont for her suggestions and correspondence which clarified many details about the divided $\be$-family at the prime $3$---we gained some insight we simply could not do without. Thank you as well to Lennart Meier for reading a draft and to Mark Behrens for his helpful comments. Thank you as well to the anonymous referees for their useful comments and suggestions. The first author was supported by NSF grant \texttt{DMS-2401918} as well as the NWO grant \texttt{VI.Vidi.193.111}. The second-named author was an associate member of the Hausdorff Center for Mathematics at the University of Bonn (\texttt{DFG GZ 2047/1}, project ID \texttt{390685813}) and is now supported by the DFG-funded research training group GRK 2240: Algebro-Geometric Methods in Algebra, Arithmetic and Topology.



\section{The height 2 connective image-of-\texorpdfstring{$J$}{J} spectrum}\label{setupsection}

The main character in this article is the $\E_\infty$-ring $\j^2$, the \emph{height $2$ connective image-of-$J$ spectrum}. Just as the classical spectrum $\j^1$ is defined as a fibre using the $p$-adic stable Adams operations on $p$-adic $K$-theory, we also need $p$-adic stable Adams operations on $\tmf$ to define $\j^2$. The following is a specialisation of \cite[Th.A \& C]{adamsontmf} by $3$-completing and taking connective covers.

\begin{theorem}\label{adamsoperationsonconnective}
    For each $k\in \Z_3^\times$, there is a map of $\E_\infty$-rings $\psi^k\colon \tmf \to \tmf$ such that the diagram of $\E_\infty$-rings
    \[\begin{tikzcd}
        {\tmf}\ar[r, "{\psi^k}"]\ar[d]  &   {\tmf}\ar[d]    \\
        {\ko}\ar[r, "{\psi^k}"]        &      {\ko}
    \end{tikzcd}\]
    commutes up to homotopy, where $\tmf\to \ko$ is the usual $q$-expansion map (see \cite[Th.A.8]{hilllawson}, \cite{laureskonetmf}, or \cite[(1.2)]{globaltate}) and $\psi^k\colon \ko\to \ko$ is the classical stable Adams operation; see \cite[\textsection5.5]{luriestheorem} for a construction of Adams operations as morphisms of $\E_\infty$-rings. Moreover, for an element $x\in \pi_\ast \tmf$, we have $\psi^k(x)=x$ if $x$ is torsion and $\psi^k(x)=k^dx$ if $x$ has degree $2d$ and lies in the chosen complement of the torsion elements $\mathfrak{F}\mathrm{ree}$ detailed in \cite[Not.3.3]{adamsontmf}.
\end{theorem}

For instance, the classes $c_4, c_6$ and $\Delta^3$ all lie in $\mathfrak{F}\mathrm{ree}$, so one computes
\[\psi^k(c_4) = k^4 c_4, \quad \psi^k(c_6) = k^6 c_6, \quad \psi^k(\Delta^3) = k^{36} \Delta^3\]
for any $k \in \Z_3^\times$; this essentially determines the $\psi^k$-action on $\pi_\ast \tmf$ by multiplicativity.\\

Our definition of $\j^2$ differs slightly from that in \cite[Df.4.9]{adamsontmf}; see \Cref{oldandnewdefinitions} for a discussion.

\begin{mydef}\label{definitionofjtwonew}
    Let $\tmf^\psi$ be the equaliser of $\psi^2$ and the identity.\footnote{The results of this article also hold if one replaces $\psi^2$ in this definition above with any $\psi^k$ for any integer $k$ not divisible by $3$. For simplicity, we stick with $k=2$.} The \emph{height 2 connective image-of-$J$ spectrum} $\j^2$ is defined using the pullback of $\E_\infty$-rings
\begin{equation}\label{definitiofclassical}\begin{tikzcd}
    {\j^2}\ar[r]\ar[d]          &   {\tau_{\geq 0}\tmf^\psi}\ar[d]   \\
    {\tau_{\leq 21}\Sph}\ar[r]   &   {\tau_{\leq 21}\tau_{\geq 0}\tmf^\psi.}
\end{tikzcd}\end{equation}
\end{mydef}

Written in this way, $\j^2$ comes with a chosen $\E_\infty$-multiplication. As the forgetful functor $\CAlg \to \Sp$ from $\E_\infty$-rings to spectra preserves limits, we can compute the homotopy groups of $\j^2$ easily using a Mayer--Vietoris sequence containing the other spectra involved in (\ref{definitiofclassical}).\\

The more fundamental definition here is $\tmf^\psi$---that is the object that implicitly appears in most of our calculations and deductions. The reason to add the $21$st truncation is to force the unit $\Sph \to \j^2$ to be more highly connected, ie, to kill some noise in the low degree homotopy groups of $\tmf^\psi$ that cannot lie in its Hurewicz image. Similar manipulations occur with the height $1$ connective image-of-$J$ spectrum $\j^1$; see \cite[\textsection4]{syntheticj} for more discussion in this direction.\\

In this section, we will discuss the homotopy groups of $\j^2$ and a construction of a synthetic spectrum $\j^2_\BP$ which implements a modified Adams--Novikov spectral sequence (ANSS) for $\j^2$.

\subsection{The homotopy groups of \texorpdfstring{$\j^2$}{j2}}\label{homotopygroupssection}
The canonical map $\j^2 \to \tau_{\leq 21} \Sph$ is $22$-connective as $\tmf^\psi \to \tau_{\leq 21} \tmf^\psi$ is. Computing $\j^2$ in low homotopy groups is therefore not a problem. In higher degrees, we use the fact that the map $\j^2 \to \tmf^\psi$ is $21$-truncated\footnote{Actually, this map is $19$-truncated.} as $\tau_{\leq 21} \Sph \to \tau_{\leq 21} \tmf^\psi$ is. In other words, $\pi_d \j^2 \to \pi_d \tmf^\psi$ is an isomorphism for $22\leq d$. One can then express $\tmf^\psi$ using the fibre sequence
\begin{equation}\label{fibredefinition}
  \tmf^\psi \to \tmf \xrightarrow{\psi^2-1} \tmf
\end{equation}
to easily compute $\pi_\ast \tmf^\psi$ and $\pi_\ast \j^2$. The only potential ambiguity is a potential extension problem.

\begin{prop}\label{extensionproblem}
    For $d$ congruent to $27$ modulo $72$, the group $\pi_d \j^2$ is $3$-torsion.
\end{prop}

One can prove this by studying $\al_1$-multiplication in this degree, but we prefer another route to set up some notation for later. Part of our proof uses a general fact about multiplication by $v_1$ in synthetic spectra modulo $3$. As we will use this fact multiple times, we state and prove it now.

\begin{lemma}\label{voneactiononreductions}
Let $X$ be a synthetic spectrum and $f\in \pi_{k,s} X$. Then we have the equality
\[\partial_0(v_1 \circ q_0(f)) = \al_1 f \in \pi_{k+3,s+1}X,\]
so $v_1\circ \overline{f}$ is of the form $\widetilde{\al_1 f}$ inside $\pi_{k+4,s}X/3$.
\end{lemma}

\begin{proof}
    The diagram of synthetic spectra
\[\begin{tikzcd}
	{\1[k,s]} & {\1/3[k,s]} & {\1/3[k-4,s]} & {\1[k-3,s-1]} \\
	X & {X/3} & {X/3[-4,0]} & {X[-3,-1]}
	\arrow["{q_0}", from=1-1, to=1-2]
	\arrow["f", from=1-1, to=2-1]
	\arrow["{v_1}", from=1-2, to=1-3]
	\arrow["{\overline{f}}", from=1-2, to=2-2]
	\arrow["{\partial_0}", from=1-3, to=1-4]
	\arrow["{\overline{f}}", from=1-3, to=2-3]
	\arrow["f", from=1-4, to=2-4]
	\arrow["{q_0}", from=2-1, to=2-2]
	\arrow["{v_1}", from=2-2, to=2-3]
	\arrow["{\partial_0}", from=2-3, to=2-4]
\end{tikzcd}\]
commutes up to homotopy, as $\overline{f}$ is defined by tensoring $f$ with $\1/3$ so the left-most and right-most squares commute, and as $\overline{f}$ is $\1/3$-linear it commutes with multiplication by $v_1$. The lower-left composite is by construction $\partial_0(v_1\circ q_0(f))$, and the upper-right composite is $\al f$ as $\partial_0(v_1)=\al$.
\end{proof}

\begin{proof}[Proof of \Cref{extensionproblem}]
    Let us prove the $d=27$ case explicitly---the other cases follow by the same arguments under $\Delta^3$-periodicity with more bookkeeping. Our extension problem comes from the exact sequence
    \[0\to (\pi_{28}\tmf)/3\simeq \F_3\{\overline{c_4^2c_6}\} \to \pi_{27} \j^2 \to \F_3\{[\al\Delta]\} \to 0\]
    which we claim splits. Indeed, if it does not split then $\dim_{\F_3}(\pi_{27}\j^2/3)=2$, and if it does split then $\dim_{\F_3}(\pi_{27}\j^2/3)=3$. An alternative calculation of $\pi_{27}\j^2/3$ uses the exact sequence
    \[\pi_{28}\tmf/3 \xrightarrow{\psi^2-1} \pi_{28}\tmf/3 \simeq \F_3\{\widetilde{\al\Delta}, \overline{c_4^2c_6}\}\xrightarrow{\partial} \pi_{27}\j^2/3 \to \pi_{27}\tmf/3\simeq \F_3\{\overline{\al\Delta}\} \xrightarrow{\psi^2-1} \pi_{27}\tmf/3,\]
    the computation of $\pi_\ast \tmf/3$ in these degrees coming from \Cref{homotopyoftmfmod3picture}. \\
    
\begin{figure}[h]
    \centering
    \begin{sseqpage}[name = homotopyoftmfmod3 , y range = {0}{7}, x range = {0}{36},axes type = frame, grid = go, yscale = 0.7, xscale = 0.32]\end{sseqpage}
    \caption{Homotopy groups of $\tmf/3$ in degrees $0$ to $36$ with Adams--Novikov filtration. Multiplication by $\overline{\al}$ and $\overline{\be}$ are drawn with black lines and multiplication by $v_1$ in orange---the region around bidegree $(28,0)$ has a more subtle $v_1$-action which is discussed in the proof of \Cref{betatwofamily}.}
    \label{homotopyoftmfmod3picture}
\end{figure}
    
    We claim that both of the maps $(\psi^2-1)$ vanish. To see this, we first notice that the action of $\psi^2$ on $\tmf$ yields $\psi^2(\overline{\al\Delta})=\overline{\al\Delta}$ and $\psi^2(\overline{c_4^2c_6})=\overline{c_4^2c_6}$, so we are left to compute $\psi^2(\widetilde{\al\Delta})$. To do this, we calculate the effect of $\psi^2$ on the class detecting $\widetilde{\al\Delta}$ on the $E_2$-page of the ANSS, which is equivalent to $\pi_{\ast,\ast}\nu\tmf/(\tau,3)$. This suffices as there are no classes of higher filtration in this stem in this ANSS. Using \Cref{voneactiononreductions} above, we have the relation $v_1\overline{\Delta} = \widetilde{\al\Delta}$ in bidegree $(28,0)$, where $v_1$ is a particular choice of lift of $\al$ coming from $\pi_{3,0}\nu\Sph/3$. We also have the relations
    \begin{equation}\label{v1andeisenstein} v_1^2 = \overline{c}_4,\qquad\qquad v_1^3=\overline{c}_6,\end{equation}
    which follow from the fact that $\overline{u}=v_1$ for $\overline{u}$ a generator of $\pi_4 \ko/3$ and the fact that for $d=8,12$ the mod $3$ reduction of the $q$-expansion map
    \[\pi_d \tmf/3 \to \pi_d \ko/3\]
    is an isomorphism; see \Cref{homotopyoftmfmod3picture}. The action of $\psi^2$ on $\pi_\ast \tmf$ shows that $\psi^2(\overline{c}_4) = \overline{c}_4$ and $\psi^2(\overline{c}_6) = \overline{c}_6$. The multiplicativity of $\psi$ then shows that $\psi^2(v_1)=v_1$ and
    \[\psi^2(\widetilde{\al\Delta})=\psi^2(v_1\overline{\Delta}) = \psi^2(v_1)\psi^2(\overline{\Delta})=v_1\overline{\Delta} = \widetilde{\al\Delta}.\]
    We then conclude that $\dim_{\F_3}\pi_{27}\j^2/3=3$, hence our sequence splits, as desired.
\end{proof}

The homotopy groups of $\j^2$ can then be read from the exact sequence on homotopy groups induced by (\ref{fibredefinition}). We display the more refined computation of the signature spectral sequence associated to the \emph{synthetic spectrum} $\j^2_\BP$, which we define shortly.

\begin{remark}[Ring structure of $\pi_\ast \j^2$]
    The multiplicative structure on $\pi_\ast\j^2$ is also easy to determine, for the most part. In particular, $\al_1$ and $\be_1$ multiplications will always, for degree reasons, stay within their ``colour'', meaning that $\al_1$- or $\be_1$-multiplication on a red dot will either be another red dot (if a black line connects them) or zero. This is simply determined from the cofibre sequence (\ref{fibredefinition}). This allows one to compute the product of blue classes with red classes. Products between red classes die, which can be seen from the multiplicativity of the Bousfield--Kan filtration of the equaliser defining $\tmf^\psi$. Examples of more interesting products, such as $\be_2\be_{6/3}$ and $\be_5\be_{6/3}$, are discussed in detail in the proof of \Cref{lowdetection}. 
\end{remark}

\begin{remark}[Comparison with definition from {\cite{adamsontmf}}]\label{oldandnewdefinitions}
    Let us denote the spectrum of \cite[Df.4.9]{adamsontmf} by $\j^2_\mathrm{old}$ for now. There is a slight difference in the definitions of $\j^2$ and $\j^2_\old$, as demonstrated by the Cartesian diagrams which define them:
    \[\begin{tikzcd}
    {\j^2}\ar[r]\ar[d]          &   {\tau_{\geq 0}\tmf^\psi}\ar[d]   \\
    {\tau_{\leq 22}\Sph}\ar[r]   &   {\tau_{\leq 22}\tau_{\geq 0}\tmf^\psi,}
\end{tikzcd}\qquad\qquad\begin{tikzcd}
    {\j^2_\old}\ar[r]\ar[d]     &   {\tau_{\geq 0}\tmf^\psi}\ar[d]   \\
    {\tau_{\leq 7}\Sph}\ar[r]   &   {\tau_{\leq 7}\tau_{\geq 0}\tmf^\psi;}
\end{tikzcd}\]
    in \cite{adamsontmf}, the second-named author forgot to include a connective cover functor above. Almost everything done here for $\j^2$ holds for $\j^2_\old$---we only choose $\j^2$ for aesthetics. Indeed, the only difference is that inside $\pi_\ast\j^2_\old$ there are nonzero classes in degrees $9$, $12$, and $20$ which do not lie in the Hurewicz image. By eliminating these classes, we obtain a spectrum $\j^2$ whose Hurewicz image looks a little nicer. For example, we can now say that an element $x\in \pi_{d}\j^2$ lies in the Hurewicz if and only if its ``$v_2^9$''-multiple does, which is not true for $\j^2_\mathrm{old}$. If one wants to study $\beta_1$-inverted phenomenæ, one should probably use a connective cover along a line of slope $1/5$.
\end{remark}

\subsection{A preferred \texorpdfstring{$\BP$}{BP}-synthetic lift of \texorpdfstring{$\j^2$}{j2}}\label{syntheticfiltrationsection}
The homotopy groups of $\j^2$ are easy to compute from \Cref{extensionproblem}, but we would also like a modified Adams--Novikov spectral sequence (ANSS) for $\j^2$ to compute its Hurewicz image. As we learned in \cite{syntheticj}, it would be much easier to define a synthetic lift of $\j^2$ whose associated spectral sequence can be determined via straightforward long exact sequences. To do this, we will use the vertical $t$-structure on $\Syn$ developed in \cite[\textsection2]{syntheticj}. In particular, taking truncations $\tau_{\leq n}^\uparrow \nu X$ in this $t$-structure roughly corresponds to cutting the ANSS for $X$ along the axis $t-s = n$ and throwing away everything to the right; see \cite[Pr.2.10 \& Rmk.2.13]{syntheticj} for more details.

\begin{mydef}\label{syntheticjdefinition}
    Let $\j^2_\BP$ be the synthetic $\E_\infty$-ring defined by pullback of synthetic $\E_\infty$-rings
    \[(\nu \tmf)^\psi = \mathrm{eq}\left(\begin{tikzcd}
        {\nu \tmf}\ar[r, "{\nu \psi^2}", shift left = 2]\ar[r, "{\id}", swap, shift right = 2]  &   {\nu\tmf}
    \end{tikzcd}\right), \qquad \begin{tikzcd}
        {\j_\BP^2}\ar[r]\ar[d]  &   {\tau^\uparrow_{\geq 0}(\nu \tmf)^\psi}\ar[d] \\
        {\tau_{\leq 21}^\uparrow \1}\ar[r]   &   {\tau_{\leq 21}^\uparrow \tau^\uparrow_{\geq 0}(\nu \tmf)^\psi.}
    \end{tikzcd}\]
\end{mydef}

Just as for $\j^2$, this definition of $\j^2_\BP$ lends itself to long exact sequence arguments. The map $\j^2_\BP \to \tau_{\leq 21}^\uparrow \1$ is $22$-connective in the vertical $t$-structure by construction, the map $\j^2_\BP \to (\nu \tmf)^\psi$ is $21$-truncated and induces an isomorphism on $\pi_{d,\ast}$ for $22\leq d$, and $(\nu\tmf)^\psi$ sits in a fibre sequence of synthetic spectra
\begin{equation}\label{definitingcofibresynthic}(\nu\tmf)^\psi \to \nu\tmf \xrightarrow{\psi^2-1} \nu\tmf.\end{equation}

\begin{prop}
    The synthetic spectrum $\j_\BP^2$ is a $\tau$-complete synthetic lift of $\j^2$.
\end{prop}

\begin{proof}
    From the facts that synthetic analogues of $\BP$-nilpotent complete spectra (e.g. connective spectra) are $\tau$-complete, see \cite[Th.4.71(1)]{filtered_van_nigtevecht}, that vertical truncations preserve $\tau$-completeness \cite[Cor.2.15]{syntheticj}, and the fact that $\tau$-complete objects are closed under limits, we see $\j_\BP^2$ is $\tau$-complete. The fact that $\tau$-inversion is exact and that the $\tau$-inversion of the diagram defining $\j^2_\BP$ is precisely (\ref{definitiofclassical}) yields an equivalence $\tau^{-1} \j^2_\BP \simeq \j^2$.
\end{proof}

Let us write $\sigma\colon \Syn \to \mathrm{Fil}(\Sp)$ for the functor which associates to each synthetic spectrum $X$ its signature spectral sequence, so the filtered spectrum associated to the $\tau$-adic tower of $X$
\[\cdots \xrightarrow{\cdot \tau}  X[0,-1] \xrightarrow{\cdot \tau} X \xrightarrow{\cdot \tau} X[0,1] \xrightarrow{\cdot \tau} \cdots.\]
The spectral sequence from the filtered spectrum $\sigma(\j_\BP^2)$, see \cite[\textsection1.2]{syntheticj}, takes the form
\begin{equation}\label{gammass}E_2= \pi_{t-s, s}(\j_\BP^2/\tau) \Longrightarrow \pi_{t-s} \j^2,\end{equation}
converging to the $\tau$-inversion of $\j^2_\BP$ given by $\j^2$. This was called the \emph{signature spectral sequence} in \cite[Df.1.5]{syntheticj}, \cite[Df.1.10]{osyn}, and \cite[Df.2.7]{smfcomputation}, but when working in $\BP$-synthetic spectra, we also call this the \emph{modified ANSS} associated to $\j^2_{\BP}$.\\

\begin{figure}[h]\begin{center}
\makebox[\textwidth]{\includegraphics[trim={3.5cm 15.5cm 1.5cm 4.2cm},clip,page = 1, scale = 1]{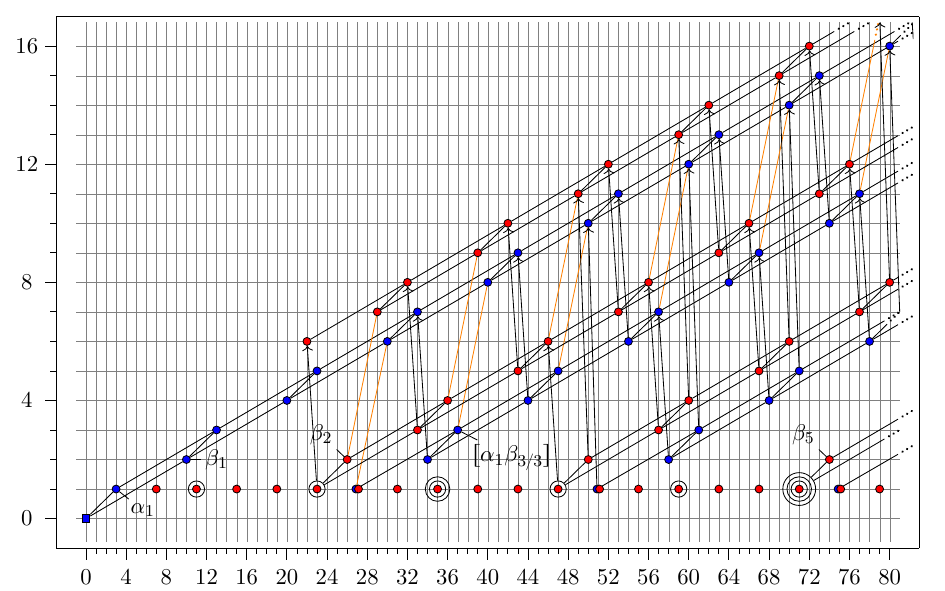}}
\caption{\label{ss1} Signature spectral sequence for $\j^2_\BP$ in stems $0\leq 80$. The classes lifted from the DSS for $\TMF$ are in blue and those in the image of the boundary map from the DSS of $\TMF$ are in red. Red and blue dots are copies of $\F_3$, the blue rectangle is a $\Z_{(3)}$, and $d$-concentric red circles are $\Z/3^d$. Orange lines indicate exotic extensions by either 3 or $\alpha_1$. We only include those red classes in filtration 1 coming from $\tmf$ given by the highest power of $\Delta$.}
\end{center}\end{figure}

\begin{figure}[h]\begin{center}
\makebox[\textwidth]{\includegraphics[trim={3.5cm 15.5cm 1.5cm 4.2cm},clip,page = 2, scale = 1]{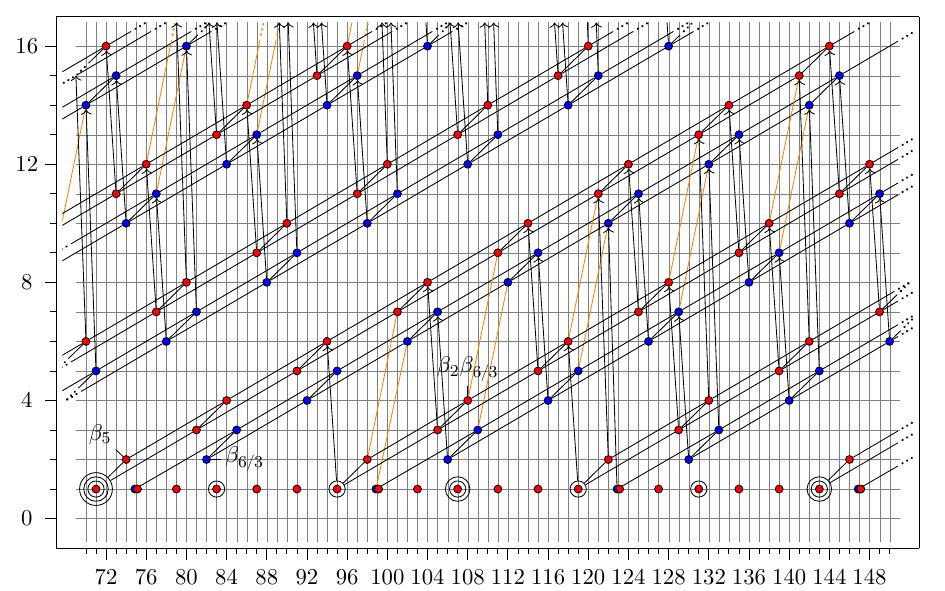}}
\caption{\label{ss2} Signature spectral sequence for $\j^2_\BP$ in stems $70\leq 150$.}
\end{center}\end{figure}

\begin{figure}[h]\begin{center}
\makebox[\textwidth]{\includegraphics[trim={3.5cm 18.5cm 1.5cm 4.2cm},clip,page = 7, scale = 1]{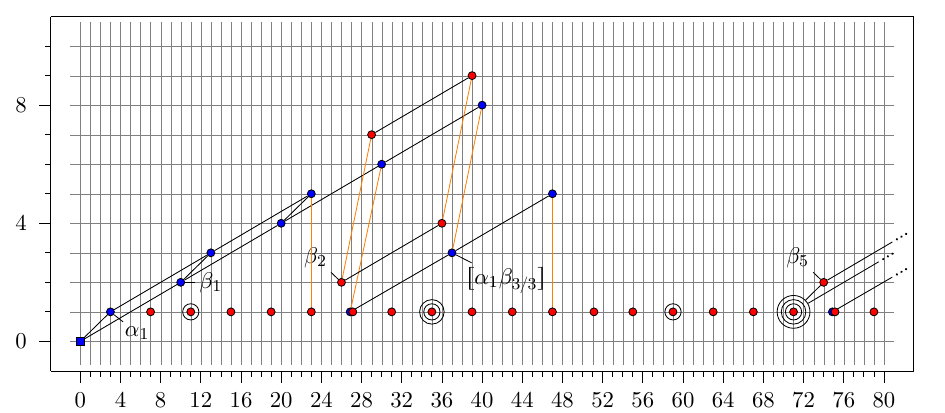}}
\caption{\label{ss3} $E_\infty$-page of signature spectral sequence for $\j^2_\BP$ in stems $0\leq 80$.}
\end{center}\end{figure}

\begin{figure}[h]\begin{center}
\makebox[\textwidth]{\includegraphics[trim={3.5cm 18.5cm 1.5cm 4.2cm},clip,page = 8, scale = 1]{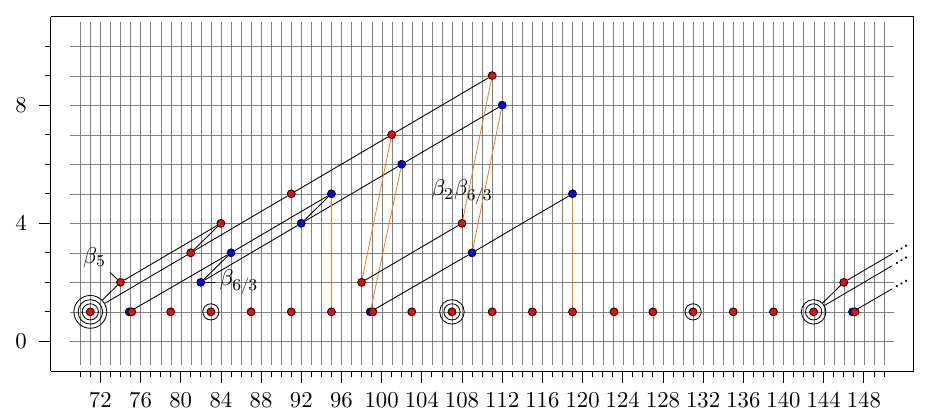}}
\caption{\label{ss4} $E_\infty$-page of signature spectral sequence for $\j^2_\BP$ in stems $70\leq 150$.}
\end{center}\end{figure}

It is easy to compute the $E_2$-page of this spectral sequence by tensoring the Cartesian diagram of \Cref{syntheticjdefinition} with $C\tau$ and computing with the associated Mayer--Vietoris exact sequences. The differentials then come from the differentials in the ANSS for $\tmf$ of \cite[\textsection7.1]{smfcomputation} and the fact that we know the homotopy groups of $\j^2$ from \Cref{homotopygroupssection}. We display the output of this spectral sequence in \Cref{ss1,ss2,ss3,ss4}.\\

What is critical about this filtration on $\pi_\ast\j^2$ is that the map $\pi_\ast\Sph \to \pi_\ast \j^2$ is one of filtered abelian groups, where $\pi_\ast\Sph$ is given its usual AN-filtration, as it comes from a map of synthetic spectra $\1\to \j_\BP^2$.

\begin{remark}[Bigraded homotopy groups of $\j^2_\BP$ and synthetic Hurewicz image]
    The bigraded homotopy groups of $\j^2_\BP$ are also easy to compute from the bigraded homotopy groups of $\nu\tmf$ or from the spectral sequence of \Cref{ss1,ss2,ss3,ss4} and the omnibus theorem \cite[Th.3.62]{filtered_van_nigtevecht}. As we will not need a full description, we will not compute them explicitly here. In the language of \cite{v232families}, one can easily compute the image of the map $\1 \to \j^2_\BP$ on bigraded homotopy groups using \Cref{maintheorem}, resulting in the so-called \emph{synthetic Hurewicz image} of $\j^2$. 
\end{remark}

\subsection{A periodic variant}
As we have done in \cite[\textsection4.4]{syntheticj}, we also want to define periodic variants of $\j^2_\BP$. The map of $\E_\infty$-rings $\psi^2\colon \TMF\to \TMF$ comes straight from ($3$-completing) \cite[Df.2.1]{heckeontmf}.

\begin{mydef}\label{periodicdefinition}
    Define the $\E_\infty$-ring $\J^2$ as the equaliser of $\psi^2$ and the identity acting on $\TMF$ and the synthetic $\E_\infty$-ring $\J^2_\BP$ as the equaliser of $\psi^2$ and the identity acting on $\nu \TMF$. We obtain natural maps $\j^2 \to \J^2$ in $\CAlg$ and $\j^2_\BP \to \J^2_\BP$ in $\CAlg(\Syn)$.
\end{mydef}

As the map $\tmf \to \TMF$ induces an isomorphism on torsion classes and an injection on torsion-free classes in nonnegative degree, a similar statement can be made for the natural map $\j^2\to \J^2$.

\begin{prop}\label{injectionintoperiodic}
    The natural map $\j^2 \to \J^2$ induces an isomorphism on classes with synthetic AN-filtration $\geq 2$ and positive stem, and is an injection otherwise. Moreover, the map $\j^2_\BP \to \J^2_\BP$ induces an injection on bigraded homotopy groups and bigraded homotopy groups modulo $\tau$. In particular, the Hurewicz image of $\j^2$ maps isomorphically onto the Hurewicz image of $\J^2$.
\end{prop}

\begin{proof}
    The first statement follows from the calculation of $\pi_\ast \j^2$ and the analogous computation of $\pi_\ast \J^2$. This implies that the Hurewicz image of $\j^2$ maps injectively into that of $\J^2$, which immediately implies the second statement too.
\end{proof}

Just as for $\j^2_\BP$, there is a signature spectral sequence associated with $\J^2_\BP$
\[E_2= \pi_{t-s, s}(\J_\BP^2/\tau) \Longrightarrow \pi_{t-s} \J^2\]
computing the homotopy groups of $\J_\BP^2$ from the homotopy groups of $\J^2_\BP/\tau$, which themselves fit in a long exact sequence with the $E_2$-page of the ANSS for $\TMF$. By \cite[Th.C]{osyn}, there is an equivalence between the ANSS for $\TMF$ and its associated descent spectral sequence, for example.

\begin{remark}[Connection to $Q(2)$]\label{qnremark}
    This periodic construction is closely related to the $Q(2)$-spectrum of Behrens \cite{ktwospheremark}. Indeed, one defines the $\E_\infty$-ring $Q(2)$ as the limit of the diagram
\begin{equation}\label{qtwodiagram}\begin{tikzcd}
    {\TMF}\ar[r, shift left = 2, "{\psi^2\times q^\ast}"]\ar[r, shift right = 2, "{\id \times p^\ast}", swap]   & {\TMF\times \TMF_0(2)}\ar[r, shift left = 3, "{\tau}"]\ar[r, "{p^\ast}"{description}]\ar[r, shift right = 3, "{\id}", swap] &   {\TMF_0(2);}
\end{tikzcd}\end{equation}
see \cite[Con.1.12]{heckeontmf} for definitions of these morphisms. By projecting onto the two $\TMF$-factors above, we obtain a map of $\E_\infty$-rings $Q(2)\to \J^2$. We predict the existence of a commutative diagram of $\E_\infty$-rings
\[\begin{tikzcd}
    {q(2)}\ar[r]\ar[d]  &   {\j^2}\ar[d]    \\
    {Q(2)}\ar[r]        &   {\J^2,}
\end{tikzcd}\]
where $q(2)$ is some connective model for $Q(2)$; we hope to use techniques similar to \cite[Th.A]{realspectra} to construct such a $q(2)$. Consequently, we see that any element that $\j^2$ detects, or equivalently by \Cref{injectionintoperiodic}, anything that $\J^2$ detects, is also detected by $q(2)$ and $Q(2)$. One can also define the synthetic $\E_\infty$-ring $Q(2)_{\BP}$ as the limit of (\ref{qtwodiagram}) after applying the synthetic analogue functor $\nu$. In future work, we will further explore these $Q(N)$-spectra from a synthetic perspective.
\end{remark}

    


\section{Detection and nondetection statements}\label{detectionsection}

With the spectrum $\j^2$ and its synthetic counterpart $\j_\BP^2$ in place, we are ready to prove \Cref{maintheorem}.

\begin{mydef}
    An element $f\in \pi_\ast\Sph$ is \emph{detected by $\j^2$} if its image under the unit map $h\colon \Sph \to \j^2$ is nonzero, i.e., if the image of $f$ lies in the Hurewicz image of $\j^2$. We say that an element $y\in\pi_d \j^2$ \emph{detects} an element $x\in \pi_d \Sph$ if $h(x)\ueq y$.
\end{mydef}

\begin{figure}[h]\begin{center}
\makebox[\textwidth]{\includegraphics[trim={3.5cm 18.5cm 2.5cm 4.2cm},clip,page = 7, scale = 1]{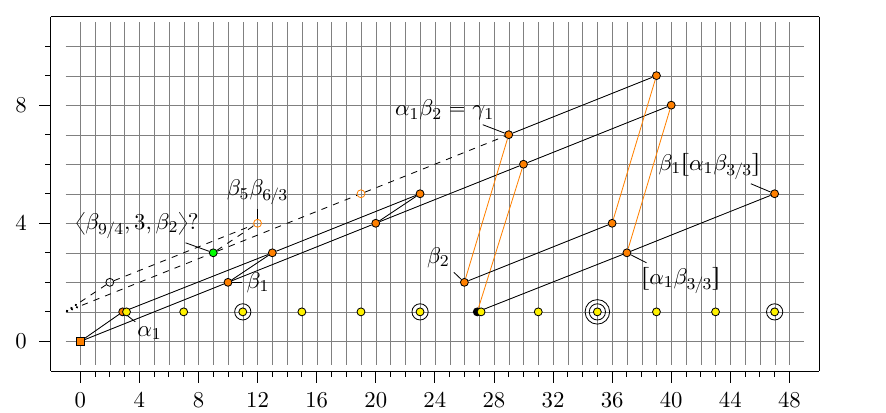}}
\caption{\label{hurewiczofjpartone} Hurewicz image of $\j^2$ in degrees $0$ to $40$ omitting those classes in filtration one which are not in the Hurewicz image and do not support multiplication by $\beta_1$. The hollow classes only appear after the $143$-stem. The green class (in the $153$-stem) lies in the Hurewicz image on $E_2$-page, but we are unsure if it lies in the Hurewicz image on homotopy groups; see \Cref{ambiguousclasses} and \Cref{153question}.}
\end{center}\end{figure}

\begin{figure}[h]\begin{center}
\makebox[\textwidth]{\includegraphics[trim={3.5cm 18.5cm 2.5cm 4.2cm},clip,page = 8, scale = 1]{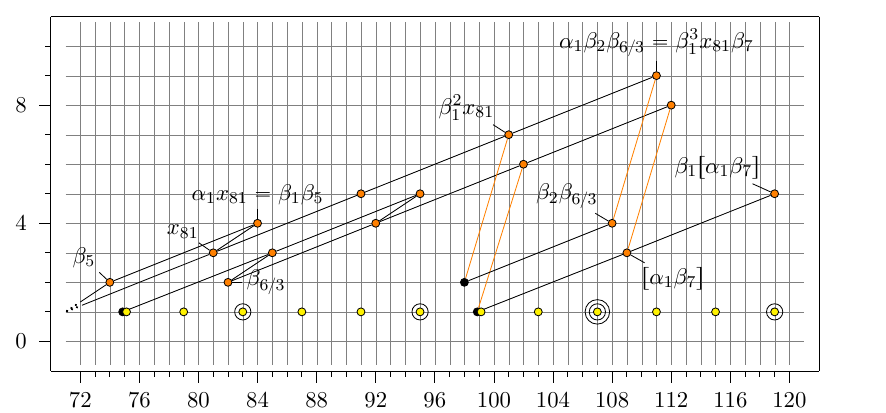}}
\caption{\label{hurewiczofjparttwo} Hurewicz image of $\j^2$ in degrees $74$ to $112$ omitting those classes in filtration one which are not in the Hurewicz image and do not support $3$- or $\be_1$-multiplication.}
\end{center}\end{figure}

The Hurewicz image, as described in \Cref{maintheorem} is displayed in orange in \Cref{hurewiczofjpartone,hurewiczofjparttwo}. Our identification of the black and orange dots in these figures proceeds as follows:
\begin{itemize}
    \item First, we recall how the orange dots corresponding to the blue part of \Cref{ss1,ss2,ss3,ss4} come from $\tmf$ and the orange dots in filtration 1 from $\j^1$; see \Cref{easyyellowdots}.
    \item Next, we show that the rest of the orange dots in degrees roughly less than $144$ lie in the Hurewicz image; see \Cref{lowdegreessection}. The pivotal proof is showing that the $\be_{2+9t}$-family is detected; see \Cref{betatwofamily}.
    \item Then orange dots are confirmed to be $144$-periodic; this is the topic of \Cref{periodicsubsection}.
    \item Finally, we show that all of the black dots do \textbf{not} lie in the Hurewicz image in \Cref{nondetectionsection}.
\end{itemize}

The green class in \Cref{hurewiczofjpartone} only definitively lies in the Hurewicz image on the $E_2$-page of the modified ANSS; this is further discussed in \Cref{ambiguousclasses} and \Cref{153question}.

\begin{remark}[Classes in filtration one]
 In \Cref{hurewiczofjpartone,hurewiczofjparttwo}, the classes in filtration $1$ are simplified for clarity. Recreating this information is straightforward. In any given positive degree $d$, the collection of elements in $\pi_d \j^2$ of filtration $1$ is isomorphic to a direct sum of $\ceil{\frac{d}{24}}$-many copies of $\pi_d \j^1$. The copy of $\pi_d \j^1$ corresponding to $\partial(c_4^ic_6^j)$ for $d=8i+12j-1$ and $0\leq j\leq 2$ is represented by the orange dot in the Hurewicz image; see \Cref{alphafamilydetection}. The rest is in the kernel of the $q$-expansion map. There are also occasionally summands supporting $\be_1$-multiplication, which correspond $\partial\Delta^{3t}$ or $\partial\al\Delta^{1+3t}$. We denote these summands with black dots in the figures above.
 \end{remark}

\begin{remark}[Connectivity of the unit $\Sph \to \j^2$]
    Unlike the height $1$ case, it is not true that $\pi_\ast \Sph \to \pi_\ast \j^2$ is surjective; this is essentially an artifact from the structure of the homotopy groups of $\tmf$. Notice that $\Sph \to \j^2$ is $22$-connective, a consequence of the fact that $\j^2\to \tau_{\leq 22} \Sph$ is $23$-connective, although it fails to be $23$-connective. Indeed, inside $\pi_{23} \j^2$ are elements not in the Hurewicz image, such as $\pm\partial 3\Delta$, which lie in the kernel of the $q$-expansion map; see \Cref{alphafamilydetection,filtrationone} for further discussion. However, $\j^2$ is still a useful approximation to $\Sph$; it detects all classes from the sphere up to $\beta_{3/2}$ in the $38$-stem as $\pi_{38}\j^2=0$.
\end{remark}


\subsection{Divided \texorpdfstring{$\alpha$}{alpha}-family and Hurewicz for \texorpdfstring{$\tmf$}{tmf}}\label{easyyellowdots}

Let us make precise the obvious detection statements from the definition of $\j^2$; this has already been summarised in \cite[Th.E]{adamsontmf}. First, we will use the canonical map $\j^2\to \tmf$ from the definition of $\j^2$.

\begin{prop}\label{hurewzicoftmf}
    The Hurewicz image of $\tmf$ is detected by $\j^2$.
\end{prop}

The elements in $\pi_\ast\tmf$ which lift through the map $\j^2\to \tmf$ are drawn in \Cref{ss1,ss2,ss3,ss4} in blue; tracing these same elements to \Cref{hurewiczofjpartone,hurewiczofjparttwo}, one finds the subset spanned by those elements in the above proposition in orange---the lower-right ``big dipper'' shape.

\begin{proof}
    The unit map $\Sph\to \tmf$ factors through $\j^2$ by construction, and the Hurewicz image of $\tmf$ is precisely the elements of (\ref{hurewiczoftmf}) by \cite[Th.6.5]{tmfthree}.
\end{proof}

Next, consider the $q$-expansion map of $\E_\infty$-rings $\tmf\to \ko$ featured in \Cref{adamsoperationsonconnective}. This sends a modular form $f$ in $\pi_\ast\tmf$ of (topological) degree $2d$ to $f(0)u^{\frac{d}{2}}$ where $u\in \pi_4\ko$ is the Bott class (recall we are implicitly working $3$-locally everywhere) and $f(0)$ is evaluation of the $q$-expansion of $f$ at $q=0$. By \Cref{adamsoperationsonconnective}, we see the Adams operation $\psi^2$ on $\tmf$ and $\ko$ agree through this $q$-expansion map. This allows us to define a map of $\E_\infty$-rings $\j^2\to \j^1$ by taking a pullback of the commutative diagram of spans of $\E_\infty$-rings
\[\begin{tikzcd}
    {\tmf^\psi}\ar[r]\ar[d] &   {\tau_{\leq 22}\tmf^\psi}\ar[d] &   {\tau_{\leq 22}\Sph}\ar[d]\ar[l]    \\
    {\ko^\psi}\ar[r] &   {\tau_{\leq 0}\ko^\psi}   &   {\tau_{\leq 0}\Sph}\ar[l]
\end{tikzcd}\]
where $\ko^\psi$ is the equaliser of $\psi^2$ and the identity on $\ko$. We also call this map $\j^2\to \j^1$ the \emph{$q$-expansion map}.

\begin{prop}\label{alphafamilydetection}
    The divided $\al$-family (\ref{hurewiczofj}) is detected by $\j^2$. In more detail, the elements $\al_{i/j}$ of order $3^j$ for $0\leq i$ and $1\leq j\leq \nu_3(i)+1$, are detected by the appropriate $3$-power multiple of $\partial(c_4^ac_6^b)$.
\end{prop}

\begin{proof}
    The chosen collection of $\partial(c_4^a c_6^b)$ spans the quotient of $\pi_\ast\j^2$ by the kernel of the $q$-expansion map in each fixed degree. It is well-known that $\j^1$ detects the family of elements $\al_{i/j}$ with $\partial_{\ko}(u^i)$, where $\partial_{\ko}\colon (\tau_{\geq 4}\ko)\to \j[1]$ is the boundary map; see \cite[Th.A]{syntheticj}, for example.
\end{proof}


\subsection{Detection of divided \texorpdfstring{$\beta$}{beta}-family elements in low degrees}\label{lowdegreessection}

To conclude that the Hurewicz image of $\j^2$ contains more than just that of $\tmf$ and $\j^1$, we will now use $\j^2_\BP$ of \Cref{syntheticjdefinition} and the associated modified ANSS (\ref{gammass}). 

\begin{prop}\label{albebe}
    For $0\leq t$, the elements $\al_1\be_1\be_{1+9t} \in \pi_{23+144t} \Sph$, $\be_1[\al_1\be_{3+9t/3}]\in \pi_{47+144t} \Sph$, and $\be_1[\al_1\be_{7+9t}] \in \pi_{119+144t} \Sph$ are detected in $\j^2$ by $\al_1\be_1^2\Delta^{6t}$, $\be_1^2[\al_1\Delta^{1+6t}]$, and $\be_1^2[\al_1\Delta^{4+6t}]$, respectively.
\end{prop}

\begin{proof}
    This follows straight from the signature spectral sequence of $\j^2_\BP$; see \Cref{ss1,ss2,ss3,ss4}. Indeed, when computing the $E_2$-page (\ref{gammass}) using the $E_2$-page of the ANSS for $\tmf$, one notices that class $\Delta$ on the $E_2$-page of ANSS for $\tmf$ is not in the kernel of $\psi^2-1$ and hence does not lift to the signature spectral sequence of $\j^2_\BP$. This means that the target of the ANSS-differential $d_5(\Delta^{1+6t}) = \pm \al\beta^2\Delta^{6t}$ for $\tmf$ is not killed in the signature spectral sequence of $\j^2_\BP$. For degree reasons, we see that $\al_1\be_1^2\Delta^{6t}$ cannot be hit by any other differentials, hence this class survives the signature spectral sequence, as desired. The same holds for the classes in degrees congruent to $47$ and $119$ modulo $144$: they are no longer hit by a differential in the modified ANSS for $\j^2$ associated with $\j^2_{\BP}$ for degree reasons.
\end{proof}

The following statement is the key to all other detection statements made here.

\begin{prop}\label{betatwofamily}
    For $0\leq t$ and $0\leq i,j\leq 1$, the elements $\al_1^i\be_1^j\be_{2+9t} \in \pi_{26+3i+10j+144t} \Sph$ are detected in $\j^2$ by $\al_1^i\be_1^j\partial\al\Delta^{1+6t}$.
\end{prop}

The proof of the above proposition assumes knowledge of the synthetic spectra and chromatic homotopy theory. In particular, for a synthetic spectrum $X$, we freely use the bigraded homotopy groups $\pi_{\ast,\ast}X/\tau$ as the $E_2$ page for a spectral sequence converging to the homotopy groups of $\tau^{-1} X$; see (\ref{gammass}) or \cite[Th.1.4]{syntheticj}. As $\j^2_\BP$ is a synthetic $\E_\infty$-ring, it comes with a unit map $\1 \to \j^2_\BP$ which induces a map between the $E_2$-page of the ANSS for $\Sph$ and the modified ANSS for $\j^2$ associated to $\j^2_\BP$. In a more chromatic direction, recall that $\be_{2+9t}$ is defined in $\pi_{26+144t} \Sph$ by \cite{marksatya}, and in particular has $E_2$-representative in the ANSS for $\Sph$ given by the image of $v_2^{2+9t} \in \pi_{32+144t,0} \1/(\tau,3,v_1)$ under the boundary maps
\[\pi_{32+144t,0} \1/(\tau,3,v_1) \xrightarrow{\partial_1} \pi_{27+144t,1} \1/(\tau,3) \xrightarrow{\partial_0} \pi_{26+144t,2} \1/\tau.\]
Notice that this expression for $\be_{2+9t}$ is \textbf{not} valid in $\1$ or $\Sph$, only modulo $\tau$. Indeed, for $t=0$, this is due to the nonzero Toda bracket $\langle \be_2, 3, \al_1\rangle = \be_1^3$ with trivial indeterminacy, by \cite[Tab.A3.4]{greenbook}, for example. By definition, this Toda bracket contains lifts of the $3$-torsion class $\beta_2$ in $\pi_{26} \Sph$ to $\pi_{27} \Sph/3$, multiplied by $v_1$, and then projected back down to $\pi_{30} \Sph$. The above Toda bracket computation shows that as $\be_1^3$ is nonzero, then any lift of $\beta_2$ to $\Sph/3$ supports multiplication by $v_1$, hence it cannot have come from $\Sph/(3,v_1)$ and therefore does not admit a further lift. The same argument works synthetically too.\\

In particular, $v_2^{9t}\circ v_2^2$ can\textbf{not} be a permanent cycle in the ANSS for $\Sph/(3,v_1)$.

\begin{proof}
    The $\al_1$- and $\be_1$-multiples follow from linearity over $\Sph$, so without loss of generality, we set $i=j=0$. We will also set $t=0$; the case for a general $1\leq t$ is the same with added bookkeeping using the fact that $\be_{2+9t}$ has $E_2$-page representative in the ANSS for $\Sph$ given by $\partial_0\partial_1(v_2^{9t}\circ v_2^2)$. Consider the diagram of abelian groups
\begin{equation}\label{bigbetwodiagram}\begin{tikzcd}
	{\pi_{32,0}\1/(\tau,3,v_1)} & {\pi_{31,1}\j^2_\BP/(\tau,3)} & {\pi_{26,2}\j_\BP^2/\tau} \\
	{\pi_{32,0}\j^2_\BP/(\tau,3,v_1)} & {\pi_{27,1}\j^2_\BP/(\tau,3)} & {\pi_{26,2}\j_\BP^2/\tau\simeq \F_3\{\partial\al\Delta\}} \\
    	& {\pi_{32,0}\j^2_\BP/(\tau,3)} & {\pi_{27,1}\j_\BP^2 / \tau \simeq \F_3\{\partial c_4^2 c_6, \al\Delta\}} \\
    {\pi_{32,0}\nu\tmf_0(2)/(\tau,3,v_1)}	& {\pi_{28,0}\j^2_\BP/(\tau,3)} & {\pi_{27,1}\j^2_\BP/\tau}
	\arrow["{\overline{\overline{h}}}"', from=1-1, to=2-1]
        \arrow[from=2-1, to=4-1]
	\arrow["{q_1}", from=3-2, to=2-1]
	\arrow["{\partial_1}", from=2-1, to=2-2]
	\arrow["{ v_1\cdot}", from=2-2, to=1-2]
	\arrow["{ v_1\cdot}", from=4-2, to=3-2]
	\arrow["{\partial_0}", from=2-2, to=2-3]
	\arrow["{q_0}", from=3-3, to=2-2]
	\arrow["{3\cdot }", from=2-3, to=1-3]
	\arrow["{3\cdot }", from=4-3, to=3-3]
\end{tikzcd}\end{equation}
where the lower-left vertical map comes from the multiplicative maps $\j^2_\BP\to \nu\tmf\to \nu\tmf_0(2)$ and the vertical zig-zags come from the cofibre sequences associated with $\j^2_\BP/(\tau,3,v_1)$ and $\j^2_\BP/(\tau,3)$. As the class $\be_2$ has an $E_2$-page representation given by $\partial_0(\partial_1(v_2^2))$, our goal is to show that 
\begin{equation}\label{v2survival}h(\be_2) = {{h}}(\partial_0(\partial_1(v_2^2))) = \partial_0(\partial_1(\overline{\overline{h}}(v_2^2))) \in \pi_{26,2}\j^2_\BP/\tau\end{equation}
is nonzero, where $h\colon \1/\tau \to \j^2_\BP/\tau$ is the unit and $\overline{\overline{h}}$ is its mod $(3,v_1)$-reduction. The latter is the $E_2$-page for the modified ANSS for $\j^2$ associated with $\j^2_\BP$, and since $\be_2$ is a permanent cycle in the ANSS for $\Sph$, it is also permanent in the modified ANSS for $\j^2$ (this is also clear in the latter for degree reasons). The desired detection now follows from the fact that $\pi_{k,\le0}\j^2_\BP/\tau=0$ for $k>0$, so that $\beta_2$ cannot be the target of a differential in the modified ANSS for $\j^2$.\\

To show (\ref{v2survival}) is nonzero, we first note that $\overline{\overline{h}}(v_2^2)$ is nonzero. Indeed, this follows from the classical fact that 
\[\tmf_0(2)\simeq \BP\langle 2\rangle \oplus \BP\langle 2\rangle [8],\]
see \cite[Pr.2.3]{hillsigmathree}, where it is credited to Hopkins--Mahowald and Behrens, so the image of $v_2^2$ in $\pi_{32,0}\nu\tmf_0(2)/(\tau,3,v_1)$ is nonzero.\\

To see that $\partial_1(\overline{\overline{h}}(v_2^2))\neq 0$, we will show that the lower $v_1$-multiplication map in (\ref{bigbetwodiagram})
\begin{equation}\label{lowvonemultiplication}
v_1\cdot\colon \F_3\{\widetilde{\al\Delta}, \widetilde{\partial c_6c_4^2}\} \simeq \pi_{28,0}\j_\BP^2/(\tau,3)\to \pi_{32,0} \j^2_\BP/(\tau,3)\simeq \F_3\{\widetilde{\partial\Delta c_4},\widetilde{\partial c_4^4}\}\end{equation}
is an isomorphism, hence $q_1$ is zero and $\partial_1$ is injective. It suffices to show (\ref{lowvonemultiplication}) is an injective map of finite dimensional vector spaces. Consider the commutative diagram of abelian groups
\[\begin{tikzcd}
    {\pi_{28,0}\j^2_\BP/(\tau,3)}\ar[r, "{v_1\cdot}"]\ar[d]  &    {\pi_{32,0}\j^2_\BP/(\tau,3)}\ar[d]  \\
    {\pi_{28,0}\nu\tmf/(\tau,3)\simeq \F_3\{\widetilde{\al\Delta}, \overline{c_6c_4^4}\}}\ar[r, "{v_1\cdot}"]    &{\pi_{32,0}\nu\tmf/(\tau,3)\simeq \F_3\{\overline{\Delta c_4}, \overline{c_4^4}\}}
\end{tikzcd}\]
from the bilinearity of tensoring with $\1/(\tau,3)$. As $\pi_{29,-1}\nu\tmf/(\tau,3)$ vanishes, the left-hand vertical map is an injection, so it suffices to show that the bottom map is an injection. Calculating on the $E_2$-page of the ANSS for $\tmf/3$ and using (\ref{v1andeisenstein}), we obtain the equalities
\[v_1 \cdot \widetilde{\al\Delta} = v_1 \widetilde{\al}\overline{\Delta} = v_1 v_1 \overline{\Delta} = v_1^2 \overline{\Delta} = \overline{\Delta c_4}, \qquad v_1 \overline{c_6 c_4^2}=\overline{c_4^4}.\]
This shows that $v_1$ is injective on $\pi_{28,0}\nu\tmf/(\tau,3)$, so that (\ref{lowvonemultiplication}) is an isomorphism.\\

The next crucial step is to identify $x = \partial_1 \overline{\overline{h}}(v_2^2)$ up to a unit. First, we notice that $x$ admits a unique lift to $\pi_{27,1}\j^2_\BP/3$ through the mod $\tau$-reduction map which is $v_1$-torsion. Indeed, to see this, first notice that the $\tau$-reduction map
\[\pi_{27,1} \j^2_\BP/3 \to \pi_{27,1} \j^2_\BP/(\tau,3) \simeq \F_3\{\widetilde{\partial\al\Delta}, \overline{\al\Delta}, \overline{\partial c_6c_4^2}\}\]
is an isomorphism as there are no classes of higher filtration; the lack of differentials in this degree also suffices. By abuse of notation, we denote elements in $\pi_{27,1} \j^2_\BP/3$ by their image modulo $\tau$. It now suffices to show that the kernel of the multiplication-by-$v_1$ map
\begin{equation}\label{uppervonemultiplucation}v_1\cdot \colon \F_3\{\widetilde{\partial\al\Delta}, \overline{\al\Delta}, \overline{\partial c_6c_4^2}\} \simeq \pi_{27,1} \j^2_\BP/3 \to \pi_{31,1} \j^2_\BP/3 \simeq \F_3\{ \overline{\partial \Delta c_4}, \overline{\partial c_4^4}, \widetilde{\be^3_1}\tau^4 \} \end{equation}
is one dimensional---we will see that this kernel is generated by $\widetilde{\partial\al\Delta}$. To begin with, by \Cref{voneactiononreductions} we obtain the equalities
\[v_1 \overline{\al\Delta} = \widetilde{\al \al\Delta} = \widetilde{\be^3_1}\tau^4\]
using the classical exotic multiplication $\al\al\Delta=\be^3\tau^4$ in $\pi_{30,2}\nu\tmf$. It was this equality that necessitates us working with $\j^2_\BP$ and not only with $\j^2_\BP/\tau$. The $\1/3$-linearity of $\overline{\partial}$ together with (\ref{v1andeisenstein}) also gets us
\[v_1 \overline{\partial c_6 c_4^2} = v_1 \overline{\partial}(\overline{c_6 c_4^2}) = \overline{\partial}(v_1 \overline{c_6 c_4^2}) = \overline{\partial}(\overline{c_4^4}) = \overline{\partial c_4^4}.\]
To finish our computation of the kernel of (\ref{uppervonemultiplucation}), we will now show that $v_1 \widetilde{\partial\al\Delta}=0$. \\

Let us start with the element $\widetilde{\be^2}\in \pi_{21} \tmf/3$. This element is uniquely determined by the fact that $\partial_0\widetilde{\be^2} = \be^2$ as $\pi_{21}\tmf =0$. Also notice that for degree reasons, we have $\widetilde{\be^2}\al=0$, where we consider $\tmf/3$ as a module over $\tmf$, and that the consequent Toda bracket $\langle \widetilde{\be^2},\al,\al\rangle$ is well-defined and has zero indeterminacy. The linearity of Toda brackets gives us the expression
\[\partial_0 \langle \widetilde{\be^2},\al,\al\rangle \subseteq \langle \partial_0 (\widetilde{\be^2}), \al, \al \rangle = \langle \be^2,\al,\al\rangle = \al\Delta \in \pi_{27} \tmf\]
which is an equality from the zero indeterminacy of the right-hand side, itself being a classical Toda bracket expression in $\tmf$; this comes from the differential $d_5(\Delta)=\pm\al\be^2$ from \cite[Pr.7.2]{smfcomputation}. In particular, $\langle \widetilde{\be^2},\al,\al\rangle$ is a fixed choice of $\widetilde{\al\Delta}$. Using this choice and the Toda juggling formula,\footnote{We need to be careful now, as $\tmf/3$ and $\Sph/3$ are no longer structured ring spectra---they are provably only $\mathbf{A}_2$; see \cite[Ex.3.3]{vig}. We are therefore considering the Toda bracket $\langle v_1, \widetilde{\be^2}, \al\rangle$ using the triangulated category $\h\Sp$ as written in \cite[Pr.2.3]{symspectravthree}. More specifically, we are considering the maps of spectra $\al\colon \Sph[3]\to \Sph$, $\widetilde{\be^2}\colon \Sph[21]\to \tmf/3$, $v_1\colon \tmf/3[4]\to \tmf/3$, with the appropriate shifts. Using the fact that $\tmf/3$ has a unital multiplication, we conclude that $v_1\circ \widetilde{\be^2} = 0$. In particular, in \emph{ibid}.\ a juggling formula is proven, which is what is used here.}
we can compute the action of $v_1$ on $\widetilde{\al\Delta}$ as
\[v_1\cdot \widetilde{\al\Delta} = v_1\cdot \langle \widetilde{\be^2}, \al, \al \rangle \ueq \langle v_1, \widetilde{\be^2}, \al \rangle \cdot {\al} \in \pi_{32}\tmf/3,\]
where $v_1\widetilde{\be^2}=0$ for degree reasons. Hence the element $v_1\widetilde{\al\Delta}$ is divisible by $\al$, which for an element in degree $\pi_{32} \tmf/3$, means that $v_1 \widetilde{\al\Delta}=0$ for degree reasons; see \Cref{homotopyoftmfmod3picture}. Let us also write $\widetilde{\al\Delta}$ for the associated cycle in $\pi_{28,0}\nu\tmf/3$, defined by applying the synthetic analogue functor $\nu$ to the map representing $\widetilde{\al\Delta}\colon \Sph[28]\to \tmf/3$. This is the correct synthetic lift as $\widetilde{\al\Delta}$ has AN-filtration $0$. Notice that this $\overline{\partial}(\widetilde{\al\Delta})$ is a choice of $\widetilde{\partial\al\Delta}$ now inside $\pi_{26,2}\j^2/3$. Indeed, this follows by chasing $\widetilde{\al\Delta}$ around the commutative diagram of abelian groups
\begin{equation}\label{littleidentificationoftilde}\begin{tikzcd}
    {\pi_{28,0}\nu\tmf/3}\ar[r, "\partial_0"]\ar[d, "\overline{\partial}"]    &   {\pi_{27,1}\nu\tmf}\ar[d, "{\partial}"]   \\
    {\pi_{27,1}\j^2_\BP/3}\ar[r, "\partial_0"]    &   {\pi_{26,2}\j^2_\BP.}
\end{tikzcd}\end{equation}
Putting everything together, we have a fixed choice of $\widetilde{\partial\al\Delta}$ and
\[v_1 \widetilde{\partial\al\Delta} = v_1 \overline{\partial}(\widetilde{\al\Delta}) = \overline{\partial}(v_1 \widetilde{\al\Delta}) = 0 \in \pi_{31,1}\j^2_\BP/3\]
from the $\1/3$-linearity of $\overline{\partial}$, as desired.\\

In summary, we know that $\partial_1(\overline{h}(v_2^2)) \ueq \widetilde{\partial\al\Delta}$ for our choice of $\widetilde{\partial\al\Delta}$ made above. By definition, $\partial_0(\widetilde{\partial \al \Delta} )= \partial\al\Delta$ is nonzero in $\pi_{26,2}\j^2_\BP/\tau$ and survives to an element of the same name in $\pi_{26}\j^2$ which detects $\be_2$ (up to a unit).
\end{proof}

The detection of the $\be_2$-family was the crucial step, from which many other families immediately follow. To see this, recall that Ravenel \cite[Table.A3.4]{greenbook} defines $x_{81} = \langle \al_1, \al_1, \be_5 \rangle$; one can check this Toda bracket has zero indeterminacy in $\pi_\ast \Sph$.

\begin{prop}\label{lowdetection}
    The elements $\al_1^i\be_2 \be_{6/3}$, $\be_1^k x_{81}$ for $0\leq i\leq 1$ and $0\leq k\leq 3$, $\al_1 x_{81}$, $\be_5$, $\be_5\be_{6/3}$, and $x_{81}\be_{6/3}$ are detected in $\j^2$ by $\al_1^i\partial(\al\be\Delta^4)$, $\be_1^{k+1} \partial(\Delta^3)$, $\al_1\be_1\partial(\Delta^3)$, $\al_1\partial(\Delta^3)$, $\al_1\be_1\partial(\Delta^6)$, and $\partial(\be^2\Delta^6)$, respectively.
\end{prop}

\begin{proof}
    To begin with, recall that both $\be_{6/3}$ and $\be_2$ are detected by $\j^2$ by \Cref{hurewzicoftmf} and \Cref{betatwofamily}, respectively. We claim that their product in $\j^2$ is the nonzero class in degree $108$ up to a unit, which we will denote by $x$ for the moment. To see this, we note that $\be_2$ and $x$ uniquely lift to classes $\al\Delta$ and $\al\be\Delta^4$ respectively under the boundary map $\partial\colon \tmf[-1] \to \j^2$, so it suffices to show that $\be_{6/3}\cdot\al\Delta = \al\be \Delta^4$ inside $\pi_{109}\tmf$ using linearity over $\Sph$. However, as $\be_{6/3}$ in $\pi_{82}\j^2$ maps to the class in $\pi_{82}\tmf$ detected by $\beta\Delta^3$, we have $\be_{6/3}\al\Delta = \al\be \Delta^4$ in $\pi_{109}\tmf$. Linearity over $\Sph$ also shows that $\j^2$ also detects $\al_1\be_2\be_{6/3}$.\\

    Using the equality $\be_2\be_{6/3} = \langle \al_1, \al_1, \be_1^2 x_{81}\rangle$ in $\pi_{108}\Sph_3$ from \cite[Table.A3.4]{greenbook}, the juggling formula for Toda brackets, and the classical Toda bracket $\be_1 = \langle \al_1, \al_1, \al_1\rangle$, we obtain the equalities
    \[\al_1 \be_2\be_{6/3} = \al_1 \langle \al_1, \al_1, \be_1^2 x_{81}\rangle = \langle \al_1, \al_1, \al_1\rangle \be_1^2 x_{81} = \be_1^3 x_{81};\]
    one may check that these Toda brackets have zero indeterminacy in $\pi_\ast\j^2$. As $\j^2$ detects $\be_1^3$, then it must also detect $x_{81}$, else the above nonzero equality would be zero. Hence $\j^2$ detects $\be_1^k x_{81}$ for $0\leq k\leq 3$. It follows that $\j^2$ detects $\al_1 x_{81}$ too. Using the Toda bracket expression $x_{81} = \langle \al_1, \al_1, \be_5\rangle$ inside $\pi_{81}\Sph$ from \cite[Table.A3.4]{greenbook} and the same tricks as above, we obtain the equalities
    \[\al_1 x_{81} = \al_1 \langle \al_1, \al_1, \be_5\rangle = \langle \al_1, \al_1, \al_1 \rangle \be_5 = \be_1 \be_5.\]
    Again, as this element and $\be_1$ are nonzero in $\pi_\ast\j^2$, then $\be_5$ must also be detected by $\j^2$. Finally, using the same argument used to show that $\be_2\be_{6/3}$ is detected in $\j^2$, one can show that $\be_5\be_{6/3}$ and $x_{81}\be_{6/3}$ are also detected in $\j^2$.
\end{proof}

Lastly, we would like to discuss the green class in bidegree $(153,3)$ in \Cref{hurewiczofjpartone}. 

\begin{remark}[$v_2^9$-periodic family from $(153,3)$]\label{ambiguousclasses}
    On the $E_2$-page of the ANSS for $\Sph$, one notices there is a class $x$ in bidegree $(153,3)$ such that $\al_1 x \ueq \be_5 \be_{6/3}$ on this $E_2$-page. In particular, $x$ must be detected on the $E_2$-page of the modified ANSS for $\j^2$ as $\j^2$ detects both $\be_5$ and $\be_{6/3}$ by \Cref{lowdetection} and \Cref{hurewzicoftmf}, respectively. If $x$ is a permanent cycle, then it would have nontrivial image in $\j^2$, but we do not know if $x$ is a permanent cycle. Moreover, we do not know if $x$ is $v_2^9$-periodic, although its image in the $E_2$-page of the modified ANSS for $\j^2$ is. 
\end{remark}

We have tried quite a lot of ideas to settle this case, but have come up empty-handed. Any suggestions towards an answer to the following question would be appreciated.

\begin{question}\label{153question}
    Is the class $x_{153,3}$ in degree $(153,3)$ in the ANSS for $\Sph$ (at the prime $3$) a permanent cycle? Is it $v_2^9$-periodic?
\end{question}

The corresponding class in the ANSS for the $\K(2)$-local sphere is a permanent cycle, but this still does not quite answer the above question.



\subsection{Detection of periodic families}\label{periodicsubsection}

The detection arguments above have calculated a lower bound for the Hurewicz image of $\j^2$ in low degrees, roughly $\leq 144$. To turn these into $144$-periodic results, we want to use the $v_2^9$-periodicity of $\tmf$. This is not so straightforward for $\j^2$ though, as $\pi_{144}\j^2=0$. Let us quickly review the height one situation as some motivation.\\

At the prime $2$ the spectrum $\j^1$ displays $8$-periodic behaviour, despite the generator of $\pi_8\j^1\simeq \F_2\{\eta \al_{4/4}\}$ being nilpotent. One often describes this through the \emph{Adams periodicity operator} $P(x)=\langle \sigma, 16, x\rangle$ applied to a $16$-torsion class $x\in \pi_d \j^1$; see \cite{adamsperiodicity}. The $16$ is necessary as $\sigma$ is $16$-torsion. If we fix a nullhomotopy of $16\sigma$ in $\j^1$, so a lift $\widetilde{\sigma}$ of $\sigma$ through the boundary map $\partial_{0^4}\colon \pi_8 \j^1/16 \to \pi_7 \j^1$, then we can further define a subset $P'(x)\subseteq P(x)$ as the set
\[P'(x) = \left\{\partial_{0^4}(\widetilde{\sigma}\cdot \widetilde{x})\ | \ \widetilde{x} \in \pi_{d+1}\j^1/16 \text{ such that }\partial_{0^4}(\widetilde{x}) = x \right\}.\]
As $\sigma = \partial(\be)$, where $\be\in \pi_8 \ko_2$ is the periodicity class and $\partial\colon (\tau_{\geq 4}\ko_2)[-1] \to \ko_2$ is the boundary map, and the reduction of $\be$ in $\pi_8 \ko/2$ is $v_1^4$, one can easily deduce that the mod $2$ reduction of $\widetilde{\sigma}$ is also $v_1^4$. From this one can show that $P'$ and also $P$ induce some kind of $v_1^4$-periodicity in the homotopy groups of $\j^1$. \\

In this subsection, we will show that the same principles are at play at height $2$ and the prime $3$, using that $\Delta^6$ is a lift of $v_2^9$ inside $\pi_{144} \tmf$. Indeed, in \Cref{v29periodicity} we will fix a lift of $\partial \Delta^6$ in $\pi_{144} \j^2/27$, in \Cref{tmfmapsagree} we will show that multiplication by this lift agrees with postcomposition by $v_2^9$ modulo $3$ and $v_1^j$, and the rest is some more calculus with, and applications of, this lift.

\begin{remark}[Periodicity operators]
    Similar to the Adams periodicity operator $P(-)=\langle \sigma,16,-\rangle$, there is also a periodicity operator acting upon the $27$-torsion inside $\j^2$ of the form $Q=\langle \partial \Delta^6, 27 , -\rangle$. Writing $Q^t$ for the $t$-fold iteration of this $Q$, when it is defined, one can check that this implements much of the $v_2^{9t}$-periodicity amongst those elements in the divided $\be$-family detected in $\j^2$:
    \[Q^t(\be_1)=\be_{1+9t},\qquad Q^t(\be_2)=\be_{2+9t},\qquad Q^t(\be_5)=\be_{5+9t}\]
    Similar to \cite[Rmk.6.7]{syntheticj}, there is some elegance to this periodicity operator, but to prove \Cref{periodicconclusions} from \Cref{betatwofamily,lowdetection} using this operator, we found ourselves repeating many of the more direct arguments found below. The reason $Q$ cannot be used as readily as $P$ or $P'$ is because $Q$ is not defined directly on the sphere. We will explore these periodicity operators in a wider context in future work.
\end{remark}

\begin{prop}\label{v29periodicity}
    There is an element $\widetilde{\partial\Delta^6}\in \pi_{144}\j^2/27$ whose image in $\pi_{144}\j^2/(3,v_1^j)$ is $v_2^9$ up to a unit for all $1\leq j\leq 4$. Moreover, this element is uniquely defined by the fact that it maps to $\partial \Delta^6$ under the boundary map $\partial_{0^3} \colon \j^2/27 \to \j^2[1]$.
\end{prop}

The number $27$ is necessary, as the group $\pi_{143}\j^2$ is $27$-torsion. We will also often think about $\widetilde{\partial\Delta^6}$ inside $\pi_{144}\j^2/3$ under the reduction map induced by the diagram of spectra
\begin{equation}\label{cofibrethreebythree}\begin{tikzcd}
    {\j^2}\ar[r, "{9\cdot}"]\ar[d, "{=}"]  &   {\j^2}\ar[r]\ar[d, "{ 3\cdot}"]  &   {\j^2/9}\ar[d]   \\
    {\j^2}\ar[r, "{ 27\cdot}"]\ar[d]  &   {\j^2}\ar[r]\ar[d]  &   {\j^2/27}\ar[d]   \\
    {0}\ar[r]  &   {\j^2/3}\ar[r, "="]  &   {\j^2/3}
\end{tikzcd}\end{equation}
whose rows and columns are cofibre sequences. The restriction $1\leq j\leq 4$ is also necessary, as the congruence up to a unit $\Delta^6 \uequiv v_2^9$ only holds modulo $(3,v_1^6)$, and the multiplication by $v_1^j$ map on $\pi_{144-4j} \tmf/3$ only has 1-dimensional cokernel for $j\leq 4$.

\begin{proof}
    First, notice that the element $\widetilde{\partial \Delta^6}$ is uniquely defined by the fact that $\partial_{0^3}(\widetilde{\partial \Delta^6})=\partial\Delta^6$ simply because the map $\partial_{0^3}$ is injective in this degree, a consequence of $\pi_{144}\j^2=0$. On the other hand, there is a class $\Delta^6$ in $\pi_{144}\tmf$ which hits $v_2^9$ up to a unit in $\pi_{144}\tmf/(3,v_1^j)$; this follows from \cite[Lm.6.2]{tmfthree}. It therefore suffices to show that the map
    \begin{equation}\label{reductiontotmfiso}\pi_{144} \j^2/(3,v_1^j)\to \pi_{144} \tmf/(3,v_1^j)\end{equation}
    is an isomorphism and sends the reduction of $\widetilde{\partial \Delta^6}$ to $\overline{\Delta^6}$ up to a unit. More specifically, we will show that (\ref{reductiontotmfiso}) is an injection, that its codomain is a $1$-dimensional $\F_3$-vector space, and that the reduction of $\widetilde{\partial\Delta^6}$ is nonzero; the reduction of $\overline{\Delta^6}$ in $\pi_{144}\tmf/(3,v_1^j)$ is necessarily nonzero as we know it is $v_2^9$ up to a unit.\\

    Contemplate the exact sequence
    \[\pi_{145} \tmf/3 \xrightarrow{q_1} \pi_{145} \tmf/(3,v_1^j) \xrightarrow{\partial_1} \pi_{144-4j} \tmf/3 \xrightarrow{ v_1^j\cdot} \pi_{144} \tmf/3.\]
    From (\ref{v1andeisenstein}) we quickly see that the right-most map is an injection, as this map is simply multiplication by $v_1$ on mod $3$ modular forms of weight $70$. Moreover, for degree reasons, we see that $\pi_{145}\tmf/3=0$. This means that $\pi_{145}\tmf/(3,v_1^j)=0$, meaning that (\ref{reductiontotmfiso}) is injective. Next, consider the exact sequence
    \begin{equation}\label{qoneissurjective}\pi_{144-4j}\tmf/3 \xrightarrow{v_1^j \cdot} \pi_{144} \tmf/3 \xrightarrow{q_1} \pi_{144} \tmf/(3,v_1^j) \xrightarrow{\partial_1} \pi_{144-4j-1} \tmf/3.\end{equation}
    The same analysis as done above shows that the cokernel of the left-most map is $\F_3\{\overline{\Delta}^6\}$, this requires $1\leq j\leq 6$, and that $\pi_{139}\tmf/3=0$. In particular, we see that $\pi_{144}\tmf/(3,v_1^j)\simeq \F_3\{\overline{\Delta}^6\}$ which shows that the codomain of (\ref{reductiontotmfiso}) is $1$-dimensional.\\
    
    Finally, to see that $q_{1^j}(\widetilde{\partial \Delta^6})$ is nonzero in $\pi_{144}\j^2/(3,v_1^j)$, it suffices to show that $\widetilde{\partial \Delta^6} \in \pi_{144} \j^2/3$ is not $v_1^j$-divisible. Well, if this class were even $v_1$-divisible, then it would be $\be_1^2$-torsion, as $\be_1^3v_1=0$ in $\pi_{34} \j^2/3 = 0$; a direct computation. However, we notice that $\be_1^3 \widetilde{\partial\Delta^6}= (\partial(\be_1^3\Delta^6))^\sim$ is nonzero, hence it cannot be $v_1$- nor $v_1^j$-divisible, and must be nonzero. This finishes the proof.
\end{proof}

It follows that powers of $\widetilde{\partial\Delta^6}$ also detect powers of $v_2^9$ in other quotients of $\j^2/3$; the following result will be useful in nondetection statements. Recall from \cite{okaringstructures}, also see \cite[Th.2.4]{tmfthree}, that $\Sph/(3,v_1^j)$ is a homotopy ring spectrum for $2\leq j$, which smashes with the $\E_\infty$-ring $\j^2$ to give the homotopy ring spectrum $\j^2/(3,v_1^j)$ for $2\leq j$.

\begin{cor}\label{congruencemodulobiggernumbers}
    Write $\widetilde{\partial \Delta^6} \in \pi_{144} \j^2/(3,v_1^j)$ for the image of the element of \Cref{v29periodicity} in $\j^2/27$. Then for $1\leq s$ and $1\leq j \leq 4s$, we have $(\widetilde{\partial \Delta^6})^{s}$ is equal to $v_2^{9s}$ in $\pi_\ast \j^2/(3,v_1^j)$ up to a unit.
\end{cor}

\begin{proof}
    This follows from \Cref{v29periodicity}; the congruence $\widetilde{(\partial \Delta^6)} \ueq v_2^9$ in $\j^2/(3,v_1^4)$ implies $\widetilde{(\partial \Delta^6)}^{9s} \ueq v_2^{9s}$ in $\j^2/(3,v_1^{4s})$.
\end{proof}

We will now utilise a little more structure on $\j^2/27$. By \cite[Th.1.2]{burklundmultiplicativestructuresmoorespectra}, $\j^2/27$ can be given an $\E_2$-structure. More precisely, we see that the natural map $\j^2/27 \to \j^2/9$ of (\ref{cofibrethreebythree}) can be refined to a map of $\E_1$-rings, also giving its fibre $\j^2/3$ the structure of a $\j^2/27$-module. The element $\widetilde{\partial \Delta^6}$ refines to a $\j^2/27$-linear self-map
\[\widetilde{\partial \Delta^6} \colon \j^2/27 \to \j^2/27[-144]\]
of the same name. Using this $\j^2/27$-linearity, we obtain the following computational lemmata.

\begin{lemma}\label{tmfmapsagree}
    For $1\leq j \leq 4$, the $\tmf$-linear maps
    \[v_2^{9}, \Delta^{6}\cdot, \widetilde{\partial\Delta^{6}}\cdot \colon  \tmf/(3,v_1^j) \to \tmf/(3,v_1^j)[-144],\]
    the first defined by smashing a $v_2^9$-self map of \cite{marksatya} or \cite{tmfthree} with $\tmf$, and the second defined by $\tmf$-scalar multiplication, and the third by $\j^2/27$-scalar multiplication, are homotopic up to a sign. Similarly, the $\j^2/27$-linear maps
    \[v_2^9, \widetilde{\partial\Delta^6}\cdot \colon \j^2/(3,v_1^j) \to \j^2/(3,v_1^j)[-144]\]
    are homotopic up to a sign.
\end{lemma}

\begin{proof}
    Consider the maps of abelian groups
    \[\left[X/(3,v_1^j), X/(3,v_1^j)[-144]\right]_{X} \xrightarrow{q_1^\ast} \left[X/3, X/(3,v_1^j)[-144]\right]_{X} \xrightarrow{q_0^\ast} \pi_{144} X/(3,v_1^j)\]
    where $X$ is either $\tmf$ or $\j^2$. We know that most of these maps agree up to a sign in the right-most group; this follows for $v_2^9$ and $\Delta^6$ on $\tmf$ from \cite[Lm.6.2]{tmfthree} and for the $X=\j^2$-cases from \Cref{v29periodicity}. The fact that $v_2^9$ and $\widetilde{\partial \Delta^6}$ agree in $\pi_{144} \j^2/(3,v_1^j)$ shows they also agree in $\pi_{144} \tmf/(3,v_1^j)$. The result follows from the fact that $q_1^\ast$ and $q_0^\ast$ are both injective, which follows from a quick computation using the defining cofibre sequences.
\end{proof}

\begin{lemma}\label{permultionvtwo}
    For $0\leq s,t$, the map $ \widetilde{\partial\Delta^6}^s\cdot \colon \pi_{75+144t} \j^2/3 \to \pi_{75+144(t+s)} \j^2/3$ induced by the $\j^2/27$-module structure on $\j^2/3$, is given by $({\partial \al\Delta^{3+6t}})^\sim \cdot \widetilde{\partial\Delta^6}^s \uequiv ({\partial \al\Delta^{3+6(s+t)}})^\sim$ up to a unit modulo $\be_1$-torsion.
\end{lemma}

\begin{proof}
    Let us set $s=1$ and $t=0$ for simplicity; the case for general $t$ follows with added bookkeeping and for general $s$ by induction. For any $0\leq k$, the class $(\partial \al\Delta^{3+6k})^\sim$ can be represented by $\overline{\partial}(({\al\Delta^{3+6k}})^\sim)$ by the same argument used to show $\overline{\partial}(\widetilde{\al\Delta}) = \widetilde{\partial \al\Delta}$; see (\ref{littleidentificationoftilde}). As the map $\overline{\partial}\colon \tmf/3[-1] \to \j^2/3$ is $\j^2/27$-linear, we see that it suffices to prove the congruence
    \[\widetilde{\al\Delta^3}\cdot \widetilde{\partial\Delta^6} \uequiv \widetilde{\al\Delta^9} \in \pi_{220} \tmf/3\]
    modulo the kernel of $\partial_0\colon \pi_{220} \tmf/3 \to \pi_{219} \tmf$, which is the subspace spanned by elements of the form $\{\overline{c_4^{2+3i} c_6 \Delta^{8-i}} \}_{0\leq i\leq 8}$. As we can rewrite $\widetilde{\al\Delta^{3k}}\ueq v_1\cdot \overline{\Delta^{3k}}$ for any $1\leq k$ by \Cref{voneactiononreductions}, we are further reduced to proving the congruence
    \[\overline{\Delta^3}\cdot \widetilde{\partial\Delta^6} \uequiv \overline{\Delta^9} \in \pi_{216} \tmf/3\]
    modulo the span of the elements $\{\overline{c_4^{3+3i}\Delta^{8-i}}\}_{0\leq i\leq 8}$ as $\overline{c_4^{2+3i}c_6 \Delta^{8-i}} = v_1 \cdot \overline{c_4^{3+3i}\Delta^{8-i}}$ by (\ref{v1andeisenstein}). We now observe that this span is also precisely the kernel of the map
    \[q_1\colon \pi_{216} \tmf/3 \to \pi_{216} \tmf/(3,v_1),\]
    as these elements are precisely the image of
    \[ v_1 \cdot\colon \pi_{212} \tmf/3 \to \pi_{216} \tmf/3.\]
    In other words, it suffices to prove the equality
    \[q_1(\overline{\Delta^9}) = \overline{\overline{\Delta^9}} ={{\Delta^6}} \cdot \overline{\overline{\Delta^3}} \ueq \widetilde{\partial\Delta^6} \cdot \overline{\overline{\Delta^3}} \in \pi_{216} \tmf/(3,v_1)\]
    up to a unit, which follows from \Cref{tmfmapsagree}.
\end{proof}

\begin{lemma}\label{permultionvtwoc}
    For $0\leq s,t$, the map $ \widetilde{\partial\Delta^6}^s\cdot \colon \pi_{82+144t} \j^2/3 \to \pi_{82+144(t+s)} \j^2/3$ induced by the $\j^2/27$-module structure on $\j^2/3$, is given by $({\partial \beta_1\Delta^{3+6t}})^\sim\cdot \widetilde{\partial\Delta^6}^s \uequiv ({\partial \beta_1 \Delta^{3+6(s+t)}})^\sim$ up to a unit modulo $\be_1$-torsion.
\end{lemma}

\begin{proof}
    This follows from \Cref{permultionvtwo} by applying the bracket operator $\langle \al_1, \al_1, - \rangle$.
\end{proof}



Using the class $\widetilde{\partial\Delta^6}$ in $\pi_{144}\j^2/27$ and the previous computation, we can show that the family $\be_{5+9t}$ is detected by $\j^2$.

\begin{prop}\label{betafivefamily}
    For $0\leq t$, the elements $\be_{5+9t}$ are detected by $\j^2$. In particular, $\be_{5+9t}$ is detected by $\al_1\partial (\Delta^{3+6t})$.
\end{prop}

\begin{proof}
    Consider the commutative diagram
    \[\begin{tikzcd}
        {\pi_{80+144t} \Sph/(3,v_1)}\ar[d, "{\overline{\overline{h}}}"]\ar[r, "{\partial_1}"]   &   {\pi_{75+144t} \Sph/3}\ar[d, "{\overline{{h}}}"]\ar[r, "{\partial_0}"]   &   {\pi_{74+144t} \Sph}\ar[d, "{{{h}}}"]    \\
        {\pi_{80+144t} \j^2/(3,v_1)}\ar[r, "{\partial_1}"]   &   {\pi_{75+144t} \j^2/3}\ar[r, "{\partial_0}"]   &   {\pi_{74+144t} \j^2\simeq \F_3.}
    \end{tikzcd}\]
    We have $\partial_0\partial_1(v_2^{5+9t})=\be_{5+9t}$ by definition, so it suffices to show that
    \[{{h}}(\be_{5+9t}) = \partial_0\partial_1(\overline{\overline{h}}(v_2^{5+9t})) = \partial_0\partial_1(v_2^{9t}\circ \overline{\overline{h}}(v_2^{5}))\]
    is nonzero; the second equality coming from the definition of the $v_2^9$ self-map on $\j^2/(3,v_1)$. We then have the equalities
    \[\partial_0\partial_1(v_2^{9t}\circ \overline{\overline{h}}(v_2^{5})) \ueq \partial_0\partial_1(\widetilde{\partial\Delta^{6}}^t \cdot \overline{\overline{h}}(v_2^{5})) \ueq \partial_0(\widetilde{\partial\Delta^{6}}^t \cdot \widetilde{\partial \al\Delta^3}) \ueq \partial_0(({\partial\al\Delta^{3+6t}})^\sim) = \partial\al\Delta^{3+6t}\]
    the first from \Cref{tmfmapsagree}, the second from the $\j^2/27$-linearity of $\partial_1$ and \Cref{lowdetection}, and the third from \Cref{permultionvtwo}.  
\end{proof}





Finally, we have to deal with the periodic family in degrees congruent to 81 modulo 144.\\

Recall that we write $x_{81} = \langle \al_1, \al_1, \beta_5 \rangle$ inside $\pi_{81} \Sph$. Ideally, and as the authors previously thought, one could define a $v_2^9$-periodic family generated by $x_{81}$ via the Toda bracket expression $\langle \al_1, \al_1, \beta_{5+9t} \rangle$ using the periodic family generated by $\beta_5$, however, we have not been able to prove that $\al_1 \beta_{5+9t}$ vanishes for $1\leq t$.

\begin{remark}
    Generally, it immediately follows from \Cref{betafivefamily} that $\langle \al_1, \al_1, \beta_{5+9t}\rangle$ is detected by $\j^2$, however, this bracket is empty if $\beta_{5+9t}$ is not $\al_1$-torsion. 
\end{remark}

In lieu of such nontrivial Toda bracket, we provide another $v_2^9$-periodic family generated by $x_{81}$.

\begin{lemma}
    The class $x_{81} \in \pi_{81} \Sph$ lifts to $\pi_{99} \Sph/(3,v_1^4)$.
\end{lemma}

\begin{proof}
    As $\pi_{81} \Sph \simeq (\Z/3)^2$, this class lifts to an element $y$ in $\pi_{82} \Sph/3$. As $\pi_{97} \Sph = \pi_{98} \Sph = 0$, then $\pi_{98} \Sph/3$ also vanishes, hence $v_1^4 y = 0$, so the desired lift to $\Sph/(3,v_1^4)$ exists.
\end{proof}

\begin{mydef}\label{df:newperiodicfamily}
    Fix a lift $z \in \pi_{99} \Sph/(3,v_1^4)$ of $x_{81}$. For $1\leq t$, write
    \[x_{81}^{(t)} = \partial_0(\partial_{1^4}( v_2^{9t} \circ z)) \in \pi_{81 + 144t} \Sph\]
    using one of the self maps $v_2^{9t} \colon \Sph/(3,v_1^4)[144] \to \Sph/(3,v_1^4)$ of \cite[Cor.4.7]{tmfthree}.
\end{mydef}

\begin{prop}\label{x81family}
    For $0\leq t$, the elements $x_{81}^{(t)}$ are detected by $\j^2$ by $\partial(\beta_1 \Delta^{3+6t})$.
\end{prop}

\begin{proof}
    The same proof as \Cref{betafivefamily} works here, so using that notation and \Cref{df:newperiodicfamily}, we have
    \[\bar{\bar{h}}(x_{81}^{(t)}) = \partial_0 \partial_{1^4}(\bar{\bar{h}}(v_2^{9t} \circ z)) = \partial_0 \partial_{1^4}(v_2^{9t} \circ \bar{\bar{h}}(z)) \ueq \partial_0 \partial_1 (\widetilde{\partial \Delta^6}^t \cdot \bar{\bar{h}}(z))\]
    \[\ueq \partial_0( (\widetilde{\partial \Delta^6}^t \cdot \widetilde{\partial(\beta_1 \Delta^3)}) \ueq \partial_0((\partial \beta_1 \Delta^{3+6t})^\sim) = \partial(\beta_1 \Delta^{3+6t}),\]
    this time appealing to \Cref{permultionvtwoc} for the second-to-last equality.
\end{proof}

This leads to a lower bound on the Hurewicz image predicted by \Cref{maintheorem}.

\begin{cor}\label{periodicconclusions}
    All of the elements described in (\ref{hurewiczofjtwoa}) and (\ref{hurewiczofjtwob}) are detected in $\j^2$.
\end{cor}

\begin{proof}
    From the multiplicative structure of $\j^2$ and \Cref{betatwofamily}, it suffices to show that the elements $\be_{5+9t}$ and $x_{81}^{(t)}$ are detected in $\j^2$, which is precisely \Cref{betafivefamily,x81family}, respectively.
\end{proof}

\subsection{Detection statements for alternative \texorpdfstring{$v_2$}{v2}-periodic families}

By \Cref{hurewzicoftmf}, we know that $\be_{1+9t}$ is detected in $\j^2$, and by \Cref{betatwofamily}, we know that $\be_{2+9t}$ is detected in $\j^2$. We can also see from the multiplicative structure of $\j^2$ that both $\be_{1}\be_{2+9t}$ and $\be_{1+9t}\be_{2}$ are both nonzero in $\pi_{36+144t}\j^2$, however, in \Cref{maintheorem}, we only referenced the first family. Many such alternative periodic families are nonzero in $\pi_\ast\j^2$ but are not included in \Cref{maintheorem}. This subsection attempts to compile an exhaustive list.

\begin{cor}\label{listofproducts}
    The families (\ref{alternativeproductstmf}) and (\ref{alternativeproductsjtwoa}) are detected by $\j^2$.
\end{cor}

\begin{remark}[Products may not be distinct]
    We do not claim that all of these elements are distinct. For example, in \cite[Th.5.3]{todabetaproducts}, Toda proves that on the $E_2$-page of the ANSS for $\Sph$, we have the equality
    \[uv\be_s\be_t = st \be_u\be_v\]
    for any $s+t=u+v$; this precise equality appears in \cite[p.624]{okashimomurabetaproducts} or \cite[Th.5.6.5]{greenbook}. This immediately shows that on this $E_2$-page we have the equalities
    \[\be_{1+9s}\be_2 = \be_1\be_{2+9s},\qquad \be_{1+9t}\be_5 = \be_1\be_{5+9t}\]
    for $0\leq s,t$. As $\j^2$ does not see anything in higher AN-filtration in these degrees, these equalities hold in $\pi_\ast \j^2$. We do not know of a reference which proves an equality similar to Toda's above involving the divided $\be$-families, such as $\be_{6+9t/3}$. We also do not know if these equalities lift from the $E_2$-page to $\pi_\ast \Sph$.
\end{remark}

\begin{proof}[Proof of \Cref{listofproducts}]
    For (\ref{alternativeproductstmf}), we will argue purely with $\tmf$; this argument could have appeared in \cite{tmfthree}, for example. If we write this product in terms of the elements in $\tmf$ which detect them, it becomes clear from the ring structure on $\pi_\ast \tmf$ that these products are nonzero. Let $0\leq a,b\leq 4$ with $a+b= 4$ and $0\leq s_a,t_b$. By \cite[Th.6.5]{tmfthree} we obtain the equalities
\[\prod \be_{1+9s_a}\be_{6+9t_b/3} \ueq \prod (\be\Delta^{6s_a}) (\be \Delta^{6t_b+3}) = \be^{4}\Delta^{\sum_a 6s_a + \sum_b (6t_b+3)} \in \pi_{\ast} \tmf\]
up to a unit. In particular, this element is nonzero in $\pi_\ast \tmf$ by inspection. The argument for (\ref{alternativeproductsjtwoa}) follows similarly, except now using the multiplicative structure on $\j^2$ and the detection statements in \Cref{maintheorem} which were proven in \Cref{easyyellowdots,lowdegreessection,periodicsubsection}. For example, for (\ref{alternativeproductsjtwoa}) we have
\[x_{81}^{(w)}\left(\prod \be_{1+9s_c} \be_{6+9t_d/3}\right) \ueq \partial (\be\Delta^{3+6w}) \left(\prod \be\Delta^{6s_c} \be \Delta^{3+6t_d}\right) = \partial (\be^4\Delta^{3+6w+\sum_c 6s_c +\sum_c (6t_d+3)})\]
which is nonzero in $\pi_\ast \j^2$ by the same arguments used to compute the multiplicative structure of $\pi_\ast j^2$ from the proof of \Cref{lowdetection}.
\end{proof}


\subsection{Nondetection statements}\label{nondetectionsection}

Having obtained a lower bound for the Hurewicz image of $\j^2$ in previous sections, we will now construct an upper bound and finish our proof of \Cref{maintheorem}. We begin in filtration 1.

\begin{prop}\label{filtrationone}
    The only elements of $\pi_\ast\j^2$ in the $1$-line in the Hurewicz image are those of \Cref{alphafamilydetection}.
\end{prop}

\begin{proof}
    As discussed in the proof of \Cref{alphafamilydetection}, the only elements in $\pi_\ast\j^2$ which can detect elements in the divided $\al$-family are those elements with nonzero image under the $q$-expansion map $\j^2\to \j^1$. In other words, precisely these $\partial(c_4^ac_6^b)$ detect the divided $\al$-family. By \cite{anoneline}, also see \cite[Th.2.2]{MRW77}, the only elements in $\pi_\ast\Sph$ with AN-filtration $\leq 1$ are in the divided $\al$-family, so there is also nothing else left for classes in $\pi_\ast\j^2$ to detect.
\end{proof}

Next, we have those classes that support multiplication by a class in the Hurewicz image to a class in the Hurewicz image.



\begin{prop}\label{filtrationtwonondetection}
    Any nonzero elements in $\pi_d\j^2$ for $d\equiv 2,98$ modulo $144$ are \textbf{not} in the Hurewicz image of $\j^2$.
\end{prop}

\begin{proof}
    By \cite{MRW77}, the only divided $\beta$-family classes that could show up on the $E_2$-page of the ANSS for the sphere in filtration $2$ and stem $d \equiv 2$ modulo $144$ are those of the form $\beta_{s3^n/j}$ with some conditions on the integers $0\leq s,n$, and $1\leq j$ with $3\nmid s$; we set the index $i$ in $\beta_{s3^n/j,i+1}$ to be $0$ here, as $\beta_{s3^n/j,i+2} = 3\beta_{s3^n/j,i+1}$ when these elements exist. We now analyse the necessary conditions on $s,n,j$ one-by-one following \cite[Th.2.6]{MRW77}. Recall that the stem of $\beta_{s3^n/j}$ is $s2^43^n - 4j - 2$:\\
    
    If $n=0$, then $j=1$ and $\beta_{s/1}$ has degree $16s-6$. As this degree is congruent to 10 modulo 16 and both 2 and 98 are congruent to 2 modulo 16, such elements cannot occur in the requisite degrees. Similarly, if $n=1$, then $j\leq 3$, and $\beta_{s3/j}$ has degree $48s-4j-2$. As this degree is congruent to $-4j-2$ modulo 48 and both 2 and 98 are congruent to 2 modulo 48, such elements cannot occur in the requisite degrees.\\

    If $2\leq n$, then $j\leq 3^n + 3^{n-1} - 1$, and $\beta_{s3^n/j}$ has degree congruent to $-4j-2$ modulo $144$. If $j=1$, we are in the wrong degree, so $2\leq j$, in which case $\j^2/(3,v_1^j)$ has a homotopy ring structure by \cite{okaringstructures}. We now have two subcases. If $3 \leq s$, then
    \[j\leq 3^n + 3^{n-1} - 1 < 3^n + 3^{n-1} \leq 3^{n-2} \cdot 12 \leq 3^{n-2} \cdot 4 s,\]
    and \Cref{congruencemodulobiggernumbers} immediately implies that the image of $\widetilde{(\partial\Delta^6)}^{s3^{n-2}}$ is congruent to $v_2^{s 3^n}$ in $\pi_\ast \j^2/(3,v_1^j)$. In particular, as $v_2^{s 3^n}$ is in the image of the reduction map $q_{1^j}$, it lies in the kernel of $\partial_{1^j}$, hence 
        \[\beta_{s3^n/j} = \partial_0 \partial_{1^j}(v_2^{s3^n})\]
    vanishes in $\pi_\ast \j^2$.\\

    On the other hand, if $1 \leq s \leq 2$, consider the following diagram of abelian groups for some $1\leq N$
        \[\begin{tikzcd}
            {\pi_{16s3^n} \j^2/(3,v_1^j)}\ar[r, "{\partial_{1^{j}}}"]\ar[d, "{v_2^{9N} \ueq \widetilde{(\partial \Delta^6)}^N}"]   &   {\pi_{16s3^n-4j-1} \j^2/3}\ar[r, "{\partial_{0}}"]\ar[d, "{\widetilde{(\partial \Delta^6)}^N}"]   &   {\pi_{16s3^n-4j-2} \j^2}    \\
            {\pi_{16s3^n+144N} \j^2/(3,v_1^j)}\ar[r, "{\partial_{1^{j}}}"]   &   {\pi_{16s3^n + 144N -4j-1} \j^2/3}\ar[r, "{\partial_{0}}"]   &   {\pi_{16s3^n + 144N -4j-2} \j^2,}
        \end{tikzcd}\]
    where the sequence commutes up to a sign by \Cref{congruencemodulobiggernumbers} for large enough $N$ and the $\j^2/27$-linearity of $\partial_{1^{j}}$. The upper composite sends $v_2^{s3^n}$ to $\beta_{s3^n/j}$, so for a contradiction, let us assume this composite is nonzero for some $j$ such that $-4j-2 \equiv 2$ or $98$ modulo $144$. In particular, $x = \partial_{1^{j}}(v_2^{s3^n})$ is nonzero and lies in the kernel of the next map in the defining exact sequence, which is multiplication by $v_1^{j}$. Just as we computed the kernel of multiplication by $v_1$ on $\pi_{27,1} \j^2_\BP/3$ (\ref{uppervonemultiplucation}) to be one dimensional spanned by $\widetilde{\partial\al\Delta}$, one computes the kernel of multiplication by $v_1^{j}$ in this degree as the one-dimensional $\F_3$-vector spaces spanned by a lift of the nonzero element $\partial(\al\Delta^M)$ for the appropriate $M$. This means that $x$ is given by this lift up to a unit, as a nonzero element in the kernel of multiplication by $v_1^{j}$, hence $\partial_0(x)$ is nonzero. Unwinding its definition, we see that $\partial_0(x) = \beta_{s3^n + 9N / j}$ is nonzero in $\j^2$. However, we are free to choose $N$ to be large enough such that the first case where $3\leq s$ kicks in, so for large enough $N$, this $\beta_{s3^n + 9N/j}$ is zero in $\j^2$, a contradiction. This finishes the proof.
\end{proof}

Combining all of the (non)detection statements for $\j^2$ yields \Cref{maintheorem}.\\

To close this section, let us detail how one can use $\j^2$ and $\tmf^\psi$ to compute the action of Adams operations on the $\tmf$-homology of $2$-cell complexes.

\begin{remark}[Adams operations on $\tmf$-homology]
    Let $x\in \pi_d \Sph$ and write $Cx$ for the $2$-cell complex $\Sph/x$. If $x$ is not detected by $\tmf$, then one easily obtains a splitting $\tmf_\ast Cx \simeq \tmf_\ast \oplus \tmf_{\ast-d-1}$. By \cite[Th.D]{realspectra}, both sides have a $\Z_3^\times$-action coming from the Adams operations of \Cref{adamsoperationsonconnective}. We claim that if $x$ happens to be detected by $\j^2$ but not $\tmf$, such as the elements in (\ref{alternativeproductsjtwoa}), then this splitting does not preserve the action on the Adams operations $\psi^k$. Let us first see this with an example. Consider $x = \be_2$, which is detected by $\j^2$ via \Cref{maintheorem} but not $\tmf$ for degree reasons. Using the notation of \Cref{definitionofjtwonew}, we have the following commutative diagram of abelian groups:
\[\begin{tikzcd}
    {0 = \tmf^{\psi^2}_{4}} & {\tmf^{\psi^2}_{30}} & {\tmf^{\psi^2}_{30}(C\beta_2)} & {\tmf^{\psi^2}_{3}} & {\tmf^{\psi^2}_{29}} \\
	0 & {\F_3\{\beta_1^3\} \simeq \tmf_{30}} & {\tmf_{30}(C\beta_2)} & {\tmf_3\simeq \F_3\{\alpha_1\}} & 0 \\
	0 & {\tmf_{30}} & {\tmf_{30}(C\beta_2)} & {\tmf_3} & 0 \\
	{\tmf^{\psi^2}_3} & {\tmf^{\psi^2}_{29}} & {\tmf^{\psi^2}_{29}(C\beta_2)} & {\tmf^{\psi^2}_{2}} & {\tmf^{\psi^2}_{28}=0}
	\arrow["{\cdot \beta_2}", from=1-1, to=1-2]
	\arrow[from=1-2, to=1-3]
	\arrow[from=1-2, to=2-2]
	\arrow["0", from=1-3, to=1-4]
	\arrow[from=1-3, to=2-3]
	\arrow["{\cdot \beta_2, \simeq}", from=1-4, to=1-5]
	\arrow[from=1-4, to=2-4]
	\arrow[from=2-1, to=2-2]
	\arrow[from=2-2, to=2-3]
	\arrow["{\psi^p-1}"', from=2-2, to=3-2]
	\arrow[from=2-3, to=2-4]
	\arrow["{\psi^2-1}"', from=2-3, to=3-3]
	\arrow[from=2-4, to=2-5]
	\arrow["{\psi^2-1}"', from=2-4, to=3-4]
	\arrow[from=3-1, to=3-2]
	\arrow[from=3-2, to=3-3]
	\arrow[from=3-2, to=4-2]
	\arrow[from=3-3, to=3-4]
	\arrow[from=3-3, to=4-3]
	\arrow[from=3-4, to=3-5]
	\arrow[from=3-4, to=4-4]
	\arrow["{\cdot \beta_2, \simeq}", from=4-1, to=4-2]
	\arrow["0", from=4-2, to=4-3]
	\arrow["\simeq", from=4-3, to=4-4]
	\arrow[from=4-4, to=4-5]
\end{tikzcd}\]

The middle column can be rewritten as the exact sequence
\[0 \to \F_3\{\be_1^3\} \to \F_3\{\be_1^3, \widetilde{\al_1}\} \xrightarrow{\psi^2-1} \F_3\{\be_1^3, \widetilde{\al_1}\} \xrightarrow{\partial} \F_3\{\partial \al_1\} \to 0\]
where $\widetilde{\al_1}$ is some lift of $\al_1$. It follows that $\psi^2(\widetilde{\al_1}) = \widetilde{\al_1} \pm \be_1^3$ for any choice of lift $\widetilde{\al_1}$ of $\al_1$. In particular, the splitting of the two centre-most rows above cannot respect Adams operations.
\end{remark}


\section{Nonvanishing results in the stable stems}\label{nonvanishingsection}

To compute the Hurewicz image of $\j^2$ above, we had to assume the computations of $\pi_{\leq 108} \Sph$ of \cite[Table.A3.4]{greenbook} in low degrees, as well as the existence of the maps $v_2^9$ on $\Sph/(3,v_1^j)$ for $1\leq j\leq 8$ of \cite{marksatya,tmfthree}. As a consequence of \Cref{maintheorem}, we immediately see that all nonzero products of elements in the Hurewicz image of $\j^2$ also detect nonzero elements in $\pi_\ast \Sph$. Not only that, but we can use some of the classes in $\pi_\ast\j^2$ which do \textbf{not} lie in the Hurewicz image to make nonvanishing statements in $\pi_\ast \Sph$ too. This is precisely what we will now explore.


\subsection{Stable homotopy groups of spheres} 

\begin{cor}\label{simplecoroarrly}
    All of the products in $\pi_\ast \Sph$ displayed in (\ref{alternativeproductstmf}) and (\ref{alternativeproductsjtwoa}) are nonzero.
\end{cor}

\begin{proof}
    By \Cref{listofproducts}, the images of all of these products in $\pi_\ast \j^2$ are nonzero.
\end{proof}

One can also leverage the fact that some classes in $\pi_\ast \j^2$ that do \textbf{not} lie in the Hurewicz image can be expressed as a Toda bracket containing only elements in the Hurewicz image. To clarify this point, let us first make a similar argument for $\tmf$, as it does not appear to be in the literature. 

\begin{cor}\label{tmftodacorollary}
    For $0\leq e,f\leq 1$ with $e+f\leq 2$ and $0\leq s_e,t_f$, the elements in $\pi_\ast\Sph$
    \[\prod \be_{1+9s_e}\be_{6+9t_f/3}\]
    all support $\al_1$-multiplication.
\end{cor}

We alternatively prove this statement as part of \Cref{simplecoroarrly} above.

\begin{proof}
    First, notice that by \cite[Th.6.5]{tmfthree} and the proof of \Cref{listofproducts}, the elements in question are all detected in $\tmf$. Let us now only detail the $\be_1\be_{6/3}$-case; the other cases are similar with added bookkeeping.\\
    
    Suppose for a contradiction that $\be_1\be_{6/3}$ was $\al_1$-torsion in $\Sph$. Then we could construct the Toda bracket $\langle \be_1\be_{6/3}, \al_1, \al_1\rangle$. From the linearity of Toda brackets with respect to maps of spectra, this would be mapped into the associated Toda bracket $\langle \be_1\be_{6/3}, \al_1, \al_1\rangle$ inside $\pi_{99} \tmf$. However, inside $\tmf$ this Toda bracket contains the associated Massey product from the $E_5$-page of the ANSS for $\tmf$ by Moss' theorem (\cite{moss,evahanamoss} or \cite[\textsection3.2]{smfcomputation}), which one calculates to be the nonzero class $\al\Delta^4$ up to a unit. The Toda bracket above also has zero indeterminacy as the group $\pi_7 \tmf$ vanishes and $\pi_{96}\tmf$ is $\al_1$-torsion. In other words, the class $\pm \al\Delta^4$ lies in the Hurewicz image of $\tmf$ as the image $\pm\langle \be_1\be_{6/3}, \al_1, \al_1\rangle$. However, by \cite[Th.6.5]{tmfthree} we know that $\al\Delta^4$ does not lie in the Hurewicz image, a contradiction.
\end{proof}

Replacing $\tmf$ with $\j^2$ in the above argument, we can prove the following nonvanishing results for elements in the divided $\be$-family.

\begin{cor}\label{todacorolalry}
    For $0\leq s,t,w$, the elements in $\pi_\ast \Sph$
    \[\be_{1+9s}\be_{1+9t}\be_{5+9w}, \qquad \be_{6+9s/3}\be_{6+9t/3}\be_{5+9w}\]
    are all nonzero.
\end{cor}

This shows the family of elements (\ref{jnondetectioncor}) are nonzero, finishing up a proof of \Cref{sphereconsequences}.

\begin{proof}
    Let us focus on the first family; the argument for the second family is similar. The product
    \[\be_{1+9s}\be_{1+9t}\be_{5+9w}\]
    is zero in $\pi_\ast\j^2$ simply for degree reasons. Assume for a contradiction that this product is also zero in $\pi_\ast\Sph$. Then we would have a well-defined Toda bracket
    \[\langle \be_{1+9s}\be_{1+9t}, \be_{5+9w}, \al_1\rangle\]
    in the homotopy groups of both $\Sph$ and $\j^2$. Moss' theorem and \Cref{ss1,ss2} again imply that this Toda bracket contains a nonzero class in $\pi_{144(s+t+w)+98}\j^2$. Moreover, the indeterminacy of this Toda bracket is zero as the group $\pi_{d}\j^2$ vanishes for $d$ congruent to $78$ modulo $144$ and is $\al_1$-torsion for $d$ congruent to $95$ modulo $144$; see \Cref{hurewiczofjparttwo}. This Toda bracket also exists in the sphere by assumption, so we see that its image in $\j^2$ lies in the Hurewicz image of $\j^2$, in contradiction with \Cref{filtrationtwonondetection}.
\end{proof}

Notice that we cannot conclude that $\be_{1+9s}\be_{6+9t/3}\be_{5+9w}$ is nonzero with these arguments, as then the associated Toda bracket picks out the class $\be_{11+9(s+t+w)}$ in $\j^2$ which \textbf{does} lie in the Hurewicz image. In \cite[Th.C]{davies_beta}, the second-named author does show that these products do not vanish, building upon the current work.\\

Combining \Cref{simplecoroarrly,todacorolalry} together yields \Cref{sphereconsequences}.

\begin{remark}[Nonexistence of periodic families]
    The computation of the Hurewicz image of $\tmf$ from \cite[Th.6.5]{tmfthree} also has an immediate implication concerning the \textbf{non}existence of the periodic families $\be_{3+9t/3}$ and $\be_{7+9t}$ in $\pi_\ast \Sph$ for $0\leq t$, recovering some cases of \cite[Ths.F \& G]{l2atprime3}. Indeed, if these families did exist, then for degree reasons, they would vanish in $\pi_\ast \tmf$. On the other hand, by \cite[Th.6.5]{tmfthree}, we know that $\al_1$ and the product of $\al_1$ with these families is nonzero in $\pi_\ast \tmf$, a contradiction. Unfortunately, we cannot see how to apply this argument to $\j^2$ to obtain any other nonexistence statements.
\end{remark}


\subsection{Homotopy groups of the mod 3 Moore spectrum}\label{ssec:Moore}

In \cite{aritashimomuromodthreemoore}, Arita--Shimomura prove various nonvanishing results of products of $\be$-family elements in $\pi_\ast \Sph/3$, the homotopy groups of the mod 3 Moore spectrum $\Sph/3$. Our analysis so far of $\j^2$ and $\Sph$ yields some generalisations of their work.\\

Recall that we write ${\be}'_{i/j} \in \pi_{16i-4j-1} \Sph/3$ for the image of $v_2^i$ under the boundary map $\partial_{1^j}\colon \Sph/(3,v_1^j) \to \Sph/3[4j+1]$, when this $v_2^i$ exists. We will also write $(\al_1\be_{i/j})'=\partial_{1^j}(\al_1v_2^i)$ when this class exists. In particular, $\partial_0(\be'_{i/j})=\be_{i/j}$ for those $i,j$ when these classes exist.

\begin{cor}\label{sillyaprime}
    For all $0\leq t$, the classes in $\pi_\ast \Sph/3$
    \[\be_{1+9t}',\qquad \be_{2+9t}',\qquad (\al_1\be_{3+9t/3})',\qquad \be_{6+9t/3}',\qquad (\al_1\be_{7+9t})'\]
    all support $\al_1$-multiplication.
\end{cor}

\begin{proof}
    Applying $\partial_0$ to these classes yields the associated nonzero classes in (\ref{alternativeproductstmf}) by the $\Sph$-linearity of $\partial_0$, so they must not vanish.
\end{proof}

The following is proven in the same way as \Cref{sillyaprime}.

\begin{cor}\label{moorespectrumnonvansihing}
    For $0\leq s,t$ such that $s\equiv 1$ modulo $9$ and $t\equiv 1,2,5$ modulo $9$, or $s\equiv 2,5$ modulo $9$ and $t\equiv 1$ modulo $9$, the element $\be_s' \be_t$ in $\pi_{16(s+t)-11} \Sph/3$ is nonzero. For all $0\leq u,v$ with $u\equiv 1,2,5$ modulo $9$, the elements $\be_{u}' \be_{6+9v/3}$, $\be_{6+9v/3}' \be_{u}$ in $\pi_{16(u+9v+6)-19} \Sph/3$ are nonzero.    
\end{cor}


One can continue to make further generalisations, for example, taking higher powers of $\be_1$; we leave such statements to the interested reader.

\begin{remark}[Connection to the literature]\label{generalisationremark}
Let us summarise the relationship between everything in this section and some other sources from the literature; further generalisations and discussions can be found in \cite{davies_beta}-

    \begin{itemize}
    \item The nonvanishing of (\ref{alternativeproductstmf}) is closely related to the first and third parts of Theorem A and the first part of Corollary B of \cite{shimbetaoneactionatthree}. In the first part of Corollary B, Shimomura proves that $\prod \be_{1+9s_a}$ is nonzero for at most $5$-fold products. We can only prove this nonvanishing for $4$-fold products, but we do allow $\be_{1+9s_a}$ to be replaced by $\be_{6+9t_b/3}$. The element $\be_{6/3}$ does occur in the third part of Theorem A in ibid., however, it seems the lack of a $v_2^9$-self map on $\Sph/(3,v_1^3)$ (which appears in \cite{tmfthree}) at the time of writing prohibited Shimomura from extending this to the family $\be_{6+9t/3}$.
    \item Similarly, the nonvanishing of (\ref{alternativeproductstmf}) is closely related to the second parts of Theorem A and Corollary B of \cite{shimbetaoneactionatthree}. Again, in the second part of Corollary B,\footnote{There seems to be a typo in the second part of Corollary B. We believe the correct statement is ``$(\prod_{i=0}^k \be_{9t_i+1})\be_{9t+2}\neq0$ if and only if $k<3$'' as more appropriately matches Theorem A in ibid.} Shimomura shows that $(\prod \be_{1+9s_i})\be_{2+9t}\neq 0$ for at most $2$-fold products. We can only handle the product $\be_{1+9s}\be_{2+9t}$, but we can also include the products $\be_{6+9s/3}\be_{2+9t}$, as well as show that these products further support $\al_1$-multiplication.
    \item The nonvanishing results in \Cref{sillyaprime,moorespectrumnonvansihing} strictly generalise Theorems C and C' from \cite{aritashimomuromodthreemoore}. Indeed, in ibid, Arita--Shimomura show that $\be_s'\be_t$ is nonzero for $s\equiv 1$ modulo $9$ and $t\equiv 1,2,5$ modulo $9$, or $s\equiv 2,5$ modulo $9$ and $t\equiv 1$ modulo $9$, and that for $s\equiv 1,2,5$ modulo $9$ the element $\be_s' \be_{6/3}$ is also nonzero. Our results above strictly generalise that. Moreover, we also recover a sliver of \cite[Th.A']{aritashimomuromodthreemoore}.\footnote{The conclusion of \Cref{sillyaprime} is comparable to Theorem A in \cite{aritashimomuromodthreemoore}. In ibid, Arita--Shimomura show that $\al_s \be_t'$ does not vanish for all $3\nmid st$. We think this is a typo. Indeed, there is no proof of Theorem A as it follows from their Theorem 3.3, which contains the extra conditions that $3\nmid t+1$ or $9| t+1$. These conditions do not hold for $t\equiv 5$ modulo $9$, which makes sense, as a simple diagram chase shows that one can choose $\beta_5'$ such that $\al_1\be_5'=0$.}
\end{itemize}
\end{remark}

\addcontentsline{toc}{section}{References}
\scriptsize
\bibliography{references} 
\bibliographystyle{alpha}

\end{document}